\documentclass[letterpaper,11pt]{article}
\usepackage{multicol, xcolor}
\usepackage{amsmath,latexsym,amsbsy,amssymb}
\usepackage{psfrag,graphicx}
\usepackage{amssymb,latexsym}
\usepackage{mathrsfs, verbatim}
\usepackage[body={15cm, 22cm}]{geometry}

\usepackage{multicol}
\usepackage{hyperref}
\usepackage{enumerate}

\numberwithin{equation}{section}

\newcommand{\qed}{\hfill$\square$\vspace{0.3truecm}}

\newcommand{\R}{\mathbb{R}}

\def\ds{\displaystyle}
\def\forall{\hbox{for all}~}

\def\argmin{\hbox{arg}\!\min}

\def\ve{\varepsilon}
\def\n{\noindent}

\def\v{\vskip 1em}

\def\bega{\begin{array}}
\def\enda{\end{array}}
\def\begi{\begin{itemize}}
\def\endi{\end{itemize}}

\def\bel{\begin{equation}\label}
\def\eeq{\end{equation}}

\def\sqr#1#2{\vbox{\hrule height .#2pt
\hbox{\vrule width .#2pt height #1pt \kern #1pt
\vrule width .#2pt}\hrule height .#2pt }}
\def\square{\sqr74}

\renewcommand{\div}{\mbox{div}\,}

\newtheorem{theorem}{Theorem}[section]
\newtheorem{corollary}[theorem]{Corollary}
\newtheorem{definition}[theorem]{Definition}

\newtheorem{remark}[theorem]{Remark}
\newtheorem{lemma}[theorem]{Lemma}
\newtheorem{proposition}{Proposition}[section]

\begin{document}
\title{\bf Metric entropy for Hamilton-Jacobi equation with uniformly directionally convex Hamiltonian}
\author{Stefano Bianchini$^{(1)}$, Prerona Dutta$^{(2)}$  and Khai T. Nguyen$^{(3)}$\\ 
\\
 {\small $^{(1)}$ SISSA, via Beirut 2, IT-34014 Trieste, ITALY }\\  
 {\small $^{(2)}$  The Ohio State University, Columbus, OH 43210, USA}\\
  {\small $^{(3)}$ North Carolina State University, Raleigh, NC 27695, USA}\\ 
 \\  {\small E-mails: ~bianchin@sissa.it,~dutta.105@osu.edu,~khai@math.ncsu.edu}
 }

\date{}

\maketitle

\begin{abstract}
The present paper  studies the BV-type regularity for viscosity solutions of the Hamilton-Jacobi equation
 \[
 u_t(t,x)+H\big(D_{\!x}  u(t,x)\big)=0,\qquad (t,x)\in (0,\infty)\times\R^d,
\]
with a coercive and uniformly directionally convex Hamiltonian $H$. More precisely, we establish a BV bound on the slope of backward characteristics $DH(D_xu(t,\cdot))$ starting at a positive time $t$. Relying on the BV bound, we quantify the metric entropy  in ${\bf W}^{1,1}_{\mathrm{loc}}\left(\R^d\right)$ for the map $S_t$ that associates to every given initial data $u_0\in{\bf Lip}\big(\R^d\big)$, the corresponding solution $S_tu_0$. Finally, a counter example is constructed to show that  both $D_xu(t,\cdot)$ and $DH(D_xu(t,\cdot))$ fail to be in $BV_{\mathrm{loc}}$ for a general strictly convex and coercive $H\in\mathcal{C}^2\left(\R^d\right)$.
\v


\end{abstract}

\section{Introduction} 
Consider a first-order Hamilton-Jacobi equation 
\begin{equation}
\label{HJ}
u_t(t,x)+H\big(D_{\!x}  u(t,x)\big)=0,\qquad (t,x)\in (0,\infty)\times\R^d,
\end{equation}
where $u:[0,+\infty)\times \R^d\to\R^d$, $D_{\!x} u =(u_{x_1},\dots, u_{x_d})$
and $H:\mathbb{R}^d\rightarrow\mathbb{R}$ is a  Hamiltonian. Due to the nonlinear dependence of the characteristic speeds on the gradient of the solution, in general a classical solution $u$ will develop singularities and the gradient $ D_{\!x}u$ will become discontinuous in finite time. To cope with this difficulty, the concept of viscosity solution was introduced by Crandall and Lions in \cite{CL} to guarantee global existence, uniqueness and stability of the Cauchy problem, under suitable assumptions on the Hamiltonian $H$. In particular, assume that
\begin{itemize}
\item [{\bf (H1)}] $H\in\mathcal{C}^1\left(\R^d\right)$ is coercive and strictly convex, i.e., $\ds \lim_{|p|\to\infty}~{H(p)\over |p|}=+\infty$ and 
\[
H(tp_1+(1-t)p_2)~<~t\cdot H(p_1)+(1-t)H(p_2),\quad t\in (0,1),~p_1, p_2\in\R^d. 
\]
\end{itemize}
The Hamilton-Jacobi equation (\ref{HJ}) generates a Hopf-Lax semigroup of viscosity solutions $\left\{S_t:\mathrm{\bf Lip}\left(\R^d\right)\to \mathrm{\bf Lip}\left(\R^d\right)\right\}_{t\geq 0}$
such that  for every  initial data $u_0\in \mathrm{\bf Lip}\left(\R^d\right)$, the corresponding unique viscosity solution of  equation (\ref{HJ}) with $u(0,x)=u_0(x)$ is computed by the Hopf-Lax representation formula 
\bel{H-L}
u(t,x)=S_t(u_0)(x)=\min_{y\in\R^d}\left\{u_0(y)+t\cdot L\left({x-y\over t}\right)\right\}
\eeq
where $L$ is the Legendre transform of $H$. In addition, if  $H$ is strongly convex then 
%, i.e., there exists a  constant $\lambda>0$ such that 
%\[
%D^2H(p)\geqslant\lambda\cdot \mathbf{I}_{d}\qquad\forall p\in\mathbb{R}^d,
%\]
%then the map $x\to u(t,x)-\ds{1\over 2\lambda t}\cdot\|x\|^2$ is concave for every $t>0$. In particular, 
$u(t,\cdot)$ is twice differentiable almost everywhere and $D_x u(t,\cdot)$ has locally bounded total variation. Thanks to  Helly's compactness theorem, the map  $S_t:\mathrm{\bf Lip}(\R^d)\to \mathrm{\bf Lip}(\R^d)$ is  compact in ${\bf W}_{\mathrm{loc}}^{1,1}\left(\R^d\right)$. A natural question arises on how to measure the degree of compactness of $S_t$. This involves  using the $\ve$-entropy, introduced by Kolmogorov and Tikhomirov in \cite{KT}:

\begin{definition}\label{DefKE}
Let $(E,\rho)$ be a metric space and $F$ be a totally bounded subset of $E$. For $\varepsilon>0$, let $\mathcal{N}_{\varepsilon}(F|E)$ 
be the minimal number of sets in a covering of $F$ by subsets of $E$ having diameter no larger than $2\varepsilon$. Then the $\varepsilon$-entropy of $F$ is defined as
\begin{equation*} 
\mathcal{H}_{\varepsilon}(F|E) ~:=~ \log_{2} \mathcal{N}_{\varepsilon}(F|E).
\end{equation*}
\end{definition}
In other words, it is the minimum number of bits needed to represent a point in a given set $F$ in the space $E$, up to an   accuracy $\varepsilon$ with respect to the metric $\rho$. Such an approach stems from a conjecture of Lax  in \cite{Lax02} for scalar conservation laws with uniformly convex fluxes. A complete answer to Lax's conjecture was provided in \cite{AON1, AON3, DLG}. This study was also extended to scalar conservation laws with  nonconvex fluxes in \cite{AON4,DNR} and to hyperbolic systems of conservation laws in \cite{AON3, AON2}. Recently, the first results on the $\ve$-entropy for sets of viscosity solutions of (\ref{HJ}) were obtained in \cite{ACN,ACN1}.  The authors proved that the minimal number of bits needed to represent a viscosity solution of (\ref{HJ})  up to an accuracy $\varepsilon$ with respect to the $\bf{W}^{1,1}$-distance is of the order $\ve^{-d}$ under the strongly convex condition on Hamiltonian $H$. There the main idea was to provide controllability results for Hamilton-Jacobi equations and a compactness result for a class of semiconcave functions. However, such a gain of BV regularity does not hold for (\ref{HJ}) with a general strictly convex Hamiltonian $H$ and the previous approach to finding $\ve$-entropy of the solution set cannot be applied. 
\medskip

\noindent In this paper, we first study the fine regularity properties of viscosity solutions to (\ref{HJ})   when  $H$ satisfies {\bf (H1)} and the following assumption of uniformly directional convexity:
\begin{itemize}
\item [{\bf (H2)}] For every constant $r>0$, it holds that
\bel{udc}
\inf_{p\neq q\in\overline{B}(0,r)}\left\langle {DH(p)-DH(q)\over |DH(p)-DH(q)|}, {p-q\over |p-q|}\right\rangle:=\lambda_r>0.
\eeq
\end{itemize}
Notice that strong convexity on  $H$ implies  {\bf (H2)} but not vice versa (e.g., $H(p)=|p|^4$). Moreover, in the scalar case ($d=1$), (\ref{udc}) holds for every $H\in\mathcal{C}^2(\R)$ with $H''>0$, without requiring strong convexity. Furthermore, we refer to Remark \ref{H-C2} which gives a sufficient condition for {\bf (H2)} in $\R^d$ with $d\geq 2$.  By the Hopf-Lax representation formula (\ref{H-L}), it is well known from \cite{CS} that the set of slopes of  backward optimal rays through $(t,x)$,
denoted by
\[
{\bf b}(t,x)=\left\{{x-y\over t}:y\in C_{t,x}\right\},\quad C_{t,x}=\argmin_{y\in\R^d}\left\{u_0(y)+t\cdot L\left({x-y\over  t}\right)\right\},
\]
reduces to a singleton ${\bf b}(t,x)=\left\{DH(D_xu(t,x))\right\}$ for almost every $(t,x)\in [0,\infty)\times \R^d$ and   can be viewed as an element in ${\bf L} ^\infty\left(\R^d\right)$. Towards the sharp estimate on $\ve$-entropy of the semigroup $S_t$, we  establish a BV bound on ${\bf b}(t,\cdot )$. 

%$M$ is a Lipschitz constant of $u_0$ then  for all $t>0$, ${\bf b}(t,\cdot)$ can be viewed as an element in ${\bf L} ^\infty\left(\R^d,\R^d\right)$ with 
%\bel{LM}
%\|{\bf b}(t,\cdot)\|_{{\bf L}^{\infty}(\R^d,\R^d)}\leq\Lambda_{M}:=\max\{|q|:L(q)\leq M|q|\}.
%\eeq

\begin{theorem}\label{BV-bound} Assume that $H$ satisfies {\bf (H1)}-{\bf (H2)}. For every $t>0$ and $u_0\in{\bf Lip}\big(\R^d\big)$ with a Lipschitz constant $M>0$,  ${\bf b}(t,\cdot)$ has locally bounded total variation and its total variation $|D_x{\bf b}(t,\cdot)|$ over an open and bounded set $\Omega\subset\R^d$ of finite perimeter is bounded by 
\[
|D_x{\bf b}(t,\cdot)|(\Omega)~\leq~{1\over \gamma_{M}}\cdot \left(\Lambda_{M}+{\mathrm{diam}(\Omega)\over t}\right)\cdot\mathcal{H}^{d-1}(\partial\Omega)+{\sqrt{d}\over t}\cdot |\Omega|
\]
for some constants $\gamma_M,\Lambda_M>0$ depending on $M$ and $H$. 
\end{theorem}
Intuitively, the uniformly directional convexity of $H$ yields a bound on the directional derivatives of ${\bf b}(t,\cdot)$ in terms of $\div_x({\bf b}(t,\cdot))$. Indeed, to prove  Theorem \ref{BV-bound}, we provide  an upper bound on  the quotient  $|D_x{\bf b}(t,\cdot)|/|\div_x {\bf b}(t,\cdot)|$ for a suitable sequence of approximate solutions, which converges uniformly and  monotonically to the given solution. In turn, the approximations of  ${\bf b}$ will  converge (in the sense that their graphs converge with respect to the Hausdorff distance). As a consequence of Theorem \ref{BV-bound}, for every $t>0$, the map $S_t:{\bf Lip}\left(\R^d\right)\to {\bf Lip}\left(\R^d\right)$ is compact in ${\bf W}^{1,1}_{{\mathrm{loc}}}\left(\R^d\right)$. 
\medskip

\noindent In the second part of the paper, we shall use the bound on the total variation of ${\bf b}(t,\cdot)$  to quantify the compactness of $S_t$ for $t>0$. More precisely, given constants $m,M,R,T>0$, consider the set of initial data
\[
\mathcal{U}_{[m,M]}:=\left\{\bar{u}\in {\bf Lip}\big(\R^d\big): |\bar{u}(0)|\leq m, \mathrm{Lip}[\bar{u}]\leq M\right\}.
\]
We establish upper and lower estimates for the $\ve$-entropy of the following solution set at time $T$ 
\[
S_{T,R}\left(\mathcal{U}_{[m,M]}\right):=\Big\{v_{\big\llcorner\square_R}~:~v\in S_T\left(\mathcal{U}_{[m,M]}\right)\Big\}
\]
with $\square_R=(-R,R)^d$ and $v_{\big\llcorner\square_R}$ denoting the restriction of $v$ on $\square_R$.
\begin{theorem} \label{second-main} Assume that $H\in\mathcal{C}^{2}\big(\R^d\big)$ and satisfies {\bf (H1)}-{\bf (H2)}. There exist constants $C_1,C_2,R_1,R_2>0$  such that for every $\ve>0$ sufficiently small, 
\[
C_1\cdot \left(\ds \Phi_M\left({\ve\over R_1}\right)\right)^{-d}\leq\mathcal{H}_{\ve}\left(S_{T,R}(\mathcal{U}_{[m,M]})\Big|{\bf W}^{1,1}\big(\square_R\big)\right)\leq C_2\cdot \left(\ds \Psi_M\left({\ve\over R_2}\right)\right)^{-d}.
\]
\end{theorem}
Here, $\Phi_M,\Psi_M$ are strictly increasing functions which depend on $H$ and will be explicitly defined in Section 4. Using Theorem \ref{BV-bound} and a result in \cite{DN}, the $\ve$-entropy of the sets of slopes of optimal rays starting at time $T$ in ${\bf L}_{\mathrm{loc}}^1\big(\R^d\big)$ is found to be of the order $\ve^{-d}$.  Thus, to achieve  the upper bound in the above theorem, we establish a quantitative relation (depending on the nonlinearity of  $H$)  between the ${\bf W}^{1,1}$-distance of  two solutions and the ${\bf L}^{1}$-distance of slopes of two corresponding  optimal rays. Finally, towards the derivation of the lower bound on $ \mathcal{H}_{\ve}\left(S_{T,R}(\mathcal{U}_{[m,M]})\Big|{\bf W}^{1,1}\big(\square_R\big)\right)$, we  study a controllability result for (\ref{HJ}). In particular,  we show that a solution to (\ref{HJ})  with a semiconvex initial condition preserves the semiconvexity on a given time interval, provided the semiconvexity constant of the initial data is sufficiently small in absolute value. 
\medskip

\noindent The remainder of this paper is organized as follows.  In Section 2, we  collect preliminary results and definitions related to semiconcave functions, BV functions and Hamilton-Jacobi equations. In Section 3, we prove the BV-type regularity for viscosity solutions. Relying on this result, in Section 4 we establish a sharp estimate on the $\ve$-entropy of the map $S_T$. Finally in Section 5, we construct a counter-example  to show that if $H\in\mathcal{C}^2(\R^d)$ satisfies {\bf (H1)} but not {\bf (H2)} then both $D_{\!x}u(t,\cdot)$ and ${\bf b}(t,\cdot)$ fail to be in $BV_{\mathrm{loc}}$ in general. 
\section{Notation and preliminaries} 
Given a positive integer $d$ and a measurable set $\Omega\subseteq\R^d$, throughout the paper we shall denote by
\begin{itemize}
\item $|\cdot|$, the Euclidean norm in $\R^d$ and 
\[
B_d(x,R)=\{y\in\R^d: |x-y|<R\}\qquad\forall R>0;
\]
\item $\langle\cdot,\cdot\rangle$, the Euclidean inner product in $\R^d$;
\item $\partial \Omega$, the boundary of $\Omega$;
\item $[x,y]$, the segment joining two points $x,y\in \mathbb{R}^d$;
\item $\# S$, the number of elements in any finite set $S$;
\item $\mathrm{Vol}(D)$, the Lebesgue measure of a measurable set $D\subset \R^d$;
\item $\omega_d :=\mathrm{Vol}(B_d(0,1))$,  the Lebesgue measure of the unit ball in $\mathbb{R}^d$;
\item $\mathbf{L}^{1}(\Omega)$, the Lebesgue space of all (equivalence classes of) summable real-valued functions on $\Omega$, equipped with the usual norm $\|\cdot\|_{\mathbf{L}^{1}(\Omega)}$ (we shall use the same symbol in case $u$ is vector-valued);
\item $\mathbf{L}^{\infty}(\Omega)$,  the space of all essentially bounded real-valued functions on $\Omega$ and $\|u\|_{\mathbf{L}^{\infty}(\Omega)}$ is the essential supremum of a function $u\in \mathbf{L}^{\infty}(\Omega)$ (we shall use the same symbol in case $u$ is vector-valued); 
\item $\mathbf{W}^{1,1}\big(\Omega)$, the  Sobolev space of functions with summable first order distributional derivatives and $\|\cdot\|_{\mathbf{W}^{1,1}(\Omega)}$ is its norm;
\item $\mathcal{C}^1(\Omega)$, the space of continuously differentiable real valued functions on $\Omega$;
\item $\mathcal{C}_c^1\big(\Omega,\R^d\big)$, the space of continuously differentiable functions $u:\Omega\to\R^d$ with a compact support;
\item $\text{\bf Lip}(\Omega)$, the space of all Lipschitz functions $f:\Omega\to\R$ and $\mathrm{Lip}[f]$ is the
Lipschitz seminorm of $f$;
\item $\mathcal{H}^{k}(E)$, the $k$-dimensional Hausdorff measure of $E\subset \R^d$;

\item For any function $f$, the function $f_{\big\llcorner\Omega}$ is the restriction of $f$ on $\Omega$;
\item ${\bf I}_d,$ the identity matrix of size $d$;
\item $\lfloor a\rfloor:=\max\{z\in\mathbb{Z}: z\leq a\}$,  the integer part $a$.
\end{itemize}

\subsection{Semiconcave and BV functions in $\R^d$} 
\subsubsection{Semiconcave functions} 
Let us recall some basic definitions and properties of semiconcave (semiconvex) functions in $\R^d$. We refer to \cite{CS} for a general introduction to
the respective theories.
\begin{definition} A continuous function $u:\Omega\to \R$ is semiconcave with a semiconcavity constant $K$ if for all $x,h\in\R^d$  with  $[x-h,x+h]\subset\Omega$, it holds that
\[
u(x+h)+u(x-h)-2u(x)\leq K\cdot |h|^2.
\] 
We say that 
\begin{itemize}
\item [-] $u$ is semiconvex  (with constant $-K$) if $-u$ is semiconcave (with constant K);
\item [-] $u$ is locally semiconcave (semiconvex) if $u$ is  semiconcave (semiconvex) in every compact set $A\subset \Omega$.
\end{itemize}
\end{definition}
We denote the distributional gradient of a semiconcave function $u$ by $Du$ and for every $x\in \Omega$ with $\Omega\subseteq\R^d$ open, we define the superdifferential and the subdifferential of $u$ at $x$ respectively by
 \begin{equation*}
 D^+u(x):=\left\{p\in\R^d: \limsup_{y\to x} {u(y)-u(x)-\langle p,y-x\rangle\over |y-x|}\leq 0\right\},
 \end{equation*}
 \begin{equation*}
 D^-u(x):=\left\{p\in\R^d: \lim\inf_{y\to x} {u(y)-u(x)-\langle p,y-x\rangle\over |y-x|}\geq 0\right\}.
 \end{equation*}
 It is clear that $D^{\pm}u(x)$ is convex and $ D^-u(x)=-D^+(-u)(x)$ for all $x\in\Omega$.
 From  \cite[Proposition 3.3.4, Proposition 3.3.10]{CS}, the superdifferential of a semiconcave function enjoys the following properties.
 
\begin{proposition}\label{con-c} Given a convex and open set $\Omega\subseteq\R^d$, let $u:\Omega\to\R$ be  semiconcave with a semiconcavity constant $K$. Then
 \begin{itemize}
 \item [(i)] The superdifferential $D^+u(x)$ is a compact, convex, nonempty set for all $x\in\Omega$. Moreover, the set-valued map $x\mapsto D^+u(x)$ is upper semicontinuous;
 \item [(ii)]  $D^+u(x)$ is a singleton if and only if $u$ is differentiable at $x$. Furthermore,  if $D^+u(x)$ is a singleton for all $x\in\Omega$, then $u\in\mathcal{C}^1(\Omega)$;
 \item [(iii)] For every $x_1,x_2\in\Omega$, it holds that 
\[
\langle p_{2}-p_{1},x_2-x_1\rangle\leq K\cdot |x_2-x_1|^2,\qquad p_i\in D^+u(x_i), i\in \{1,2\}.
\]
 \end{itemize}
 \end{proposition}
From (ii) if $u$ is both locally semiconcave and locally semiconvex then $u$ is in $\mathcal{C}^1(\Omega)$ as shown in \cite[Corollary 3.3.8]{CS}. This is crucial to prove further regularity for viscosity solutions in Proposition \ref{conv-conx} which allows us to construct a backward smooth solution of  (\ref{HJ}). To complete this part, for every constant  $r,K>0$, let us define the set
\bel{SCrk}
\mathcal{SC}_{[r,K]}:=\left\{v\in {\bf Lip}(\R^d): \mathrm{Lip}[v]\leq r~\mathrm{and}~ v~\mathrm{is~semiconcave~with~constant}~K\right\}.
\eeq
From the  proof of \cite[Proposition 10]{ACN}, one can easily obtain a lower bound on the $\ve$-entropy for the set $\left\{D_{\!x}v_{\big\llcorner\square_R}:v\in\mathcal{SC}_{[r,K]}\right\}$ in ${\bf L}^1\left(\square_R\right)$ which  will be used to establish a lower estimate on the $\ve$-entropy of a set of viscosity solutions in subsection \ref{LUP}.
\medskip

\begin{corollary}\label{low} Given any $r,R,K>0$, for every $\ve>0$ sufficiently small, 
%\[
%0~<~\ve~\leq~\min \{r,K\}\cdot {\omega_d\cdot R^d\over (d+1)2^{d+8}},
%\]
there exists a subset $\mathcal{G}^{R}_{[r,K]}$ of $\mathcal{SC}_{[r,K]}$ such that
\[
\# \mathcal{G}^{R}_{[r,K]}~\geq~2^{ \beta_{[R,K]}\cdot \ve^{-d}},\qquad \beta_{[R,K]}~=~{1\over 3^{d}2^{d^2+4d+3}\ln 2}\cdot \left(K\omega_d R^{d+1}\over(d+1)\right)^d
\]
and
\[
\left\|D_{\!x}v_{\big\llcorner\square_R}-D_{\!x}w_{\big\llcorner\square_R}\right\|_{{\bf L}^1\big(\square_R\big)}~\geq~2\ve\qquad\forall v\neq w\in \mathcal{G}^{R}_{[r,K]}.
\]
\end{corollary}
\subsubsection{Functions of bounded total variation}

Let us now introduce the concept of functions of bounded variation. We refer to \cite{ANP} for a comprehensive analysis on this topic.
\begin{definition} The function $u\in {\bf L}^1(\Omega)$ is a function of bounded variation on $\Omega\subseteq\R^d$ and said to be in $BV(\Omega,\R^m)$, if the distributional derivative of $u$, denoted by $Du$, is an  $m\times d$ matrix of finite measures $D_iu^{\alpha}$ in $\Omega$ satisfying
\[
\sum_{\alpha=1}^m\int_{\Omega}~u^{\alpha}\div\varphi^{\alpha}~dx~=~-\sum_{\alpha=1}^m\sum_{i=1}^d\int_{\Omega}\varphi_i^{\alpha} dD_iu^{\alpha}\quad\forall \varphi\in\left[\mathcal{C}_c^1(\Omega,\R^d)\right]^m.
\]
We denote by $|Du|(\Omega)$ the total variation of $u$ over $\Omega$, i.e.,
\[
|Du|(\Omega)~=~\sup\left\{\sum_{\alpha=1}^{m}\int_{\Omega}u^{\alpha}\div^{\alpha} \varphi:\varphi\in\left[\mathcal{C}_c^1(\Omega,\R^d)\right]^m, \|\varphi\|_{\infty}\leq 1\right\}.
\]
\end{definition}
We recall a Poincar\'e-type inequality for bounded total variation functions on convex domains that will be used in the paper. This result is based on \cite[Theorem 3.2]{AD}. 
\medskip

%and on \cite[Proposition 3.2.1, Theorem 3.44]{ANP}.
\begin{theorem}(Poincar\'e inequality)\label{PC} Let $\Omega\subset \R^d$ be an open,  bounded, convex set with Lipschitz boundary.
%\item For any $u\in {\bf W}^{1,1}(\Omega)$, it holds 
%\[
%\int_{\Omega}|u|dx~\leq~|\Omega|^{1\over d}\cdot \int_{\Omega}|D u|dx.
%\]
For any  $u\in BV(\Omega,\R)$, it holds that
\[
\int_{\Omega} \big|u(x)-u_{\Omega}\big|~dx~\leq~{\mathrm{diam}(\Omega)\over 2}\cdot |Du|(\Omega)
\]
where 
\[
u_{\Omega}~=~{1\over \mathrm{Vol}(\Omega)}\cdot \int_{\Omega}u(x)~dx
\]
is the mean value of $u$ over $\Omega$.
\end{theorem}
To complete this subsection, let us recall a result on the metric entropy for a class of functions with bounded total variation which will be used  in subsection \ref{UPP}. For every $R,M,V>0$,  we consider a class of uniformly bounded total variation functions on $\square_R$
\bel{F}
\mathcal{F}_{[R,M,V]}~=~\left\{f:\square_R\to\R^d:\|f\|_{{\bf L}^{\infty}\big(\square_R\big)}\leq M, |Df|\big(\square_R)\leq V \right\}.
\eeq
By a slight modification in the proof of \cite[Theorem 1]{DN}, one can obtain the following upper bound on the $\ve$-entropy of $\mathcal{F}_{[R,M,V]}$ in ${\bf L}^1\left(\square_R\right)$.

\begin{corollary}\label{e-p} For every $0<\ve<\min\left\{{6RV^2\over 3V+2M}, 2V \left({RV\over M}\right)^{1\over d}\right\}$, it holds that
\bel{m-est}
\mathcal{H}_{\ve}\left(\mathcal{F}_{[R,M,V]}~\Big|~{\bf L}^1(\square_R)\right)~\leq~48\sqrt{d}\cdot \left({6d\sqrt{d}RV\over \ve}\right)^d.
\eeq
\end{corollary}
{\bf Proof.} By the definition of $\ve$-entropy, we have 
\bel{F1}
\mathcal{H}_{\ve}\left(\mathcal{F}_{[R,M,V]}~\Big|~{\bf L}^1(\square_R)\right)~\leq~d\cdot \mathcal{H}_{{\ve\over d}}\left(\mathcal{F}^1_{[R,M,V]}~\Big|~{\bf L}^1(\square_R)\right)
\eeq
with 
\[
\mathcal{F}^1_{[R,M,V]}~=~\left\{f:\square_R\to\R:\|f\|_{{\bf L}^{\infty}\big(\square_R\big)}\leq M, |Df|\big(\square_R)\leq V \right\}.
\]
Consider a class of real-valued bounded total variation functions
\[
\mathcal{B}_{[R,M,V]}~=~\left\{g:[0,R]\to [0,M]: |Dg|([0,R])\leq V \right\}.
\]
From \cite[Lemma 2.3]{DN}, for every $0<\ve< \ds{RV\over 3}$, one has
\[
\mathcal{N}_{\ve}\left(\mathcal{B}_{[R,{9\over 8}V,V]}\Big| {\bf L}^1([0,R])\right)~\leq~2^{{17RV\over \ve}}
\]
and this implies that 
\[
\mathcal{N}_{\ve}\left(\mathcal{B}_{[R,M,V]}\Big| {\bf L}^1([0,R])\right)~\leq~{8M\over V}\cdot \mathcal{N}_{\ve}\left(\mathcal{B}_{[R,{9\over 8}V,V]}\Big| {\bf L}^1([0,R])\right)~\leq~{8M\over V}\cdot 2^{{17RV\over \ve}}.
\]
In particular, for every $0<\ve<\ds{RV^2\over 3V+M}$ such that $\ds{8M\over V}\leq 2^{RV\over \ve}$, it holds that
\[
\mathcal{H}_{\ve}\left(\mathcal{B}_{[R,M,V]}\Big| {\bf L}^1([0,R])\right)~=~\log_2\left(\mathcal{N}_{\ve}\left(\mathcal{B}_{[R,M,V]}\Big| {\bf L}^1([0,R])\right)\right)~\leq~{18RV\over \ve}.
\]
Using the above estimate, one can follow the same argument as in the proof of  \cite[Theorem 3.1]{DN} to obtain that for every $0<\ve<\min\left\{\ds{6RV^2\over 3V+2M}, 2V \left({RV\over M}\right)^{1\over d}\right\}$, it holds that 
\[
\mathcal{H}_{\ve}\left(\mathcal{F}^1_{[R,M,V]}~\Big|~{\bf L}^1(\square_R)\right)~\leq~{48\over \sqrt{d}}\cdot \left({6\sqrt{d}RV\over \ve}\right)^d
\]
and then (\ref{F1}) yields (\ref{m-est}).
\qed
\subsection{Semigroup of Hamilton-Jacobi equation}
Consider the Hamilton-Jacobi equation~\eqref{HJ} under the assumptions {\bf (H1)}-{\bf (H2)}. Without loss of generality, we shall assume that the Hamiltonian satisfies further conditions
\bel{Hs}
H(0)=0\qquad\mathrm{and}\qquad DH(0)=0,
\eeq
otherwise the transformations $x\mapsto x+tD H(0)$, $u(t,\cdot )\mapsto u(t,x)+t\cdot H(0)$ and $H(p)\mapsto H(p)-\langle D H(0),p\rangle$  reduce the general case to this one. Before recalling the concept of viscosity solution to (\ref{HJ}), let us give a sufficient condition for the assumption  {\bf (H2)}.
\begin{remark}\label{H-C2}  Let $H$ be in $\mathcal{C}^2(\R^d)$. Assume that there exists a constant $\lambda>0$ such that 
\bel{udc2}
D^2H(p)~=~ |D^2H(p)|\cdot A(p),\qquad A(p)~\geq~\lambda\cdot \mathbf{I}_d,
\eeq
with $A(p)$  being a $d\times d$ matrix and $|D^2H(p)|$ denoting the matrix norm of $D^2H(p)$. Then $H$ satisfies {\bf (H2)}. 
\end{remark}
{\bf Proof.}
For any $p\neq q\in\R^d$, by mean value theorem, it holds that
\begin{equation*}
\begin{split}
D H(p)-D H(q) &=~\int_{0}^1D^2H(tp+(1-t)q)\cdot (p- q)dt \\[2mm]
&=~\bigg[ \int_0^1 A(t p+(1-t) q  ) \big|D^2 H (tp+(1-t)q) \big|dt \bigg] \cdot (p- q) 
%\\[2mm]
%&\doteq~\tilde A(p,q) \cdot (p - q).
\end{split}
\end{equation*}
and 
\[
|D H(p)-D H(q)|~\leq~|p-q|\cdot \int_{0}^1|D^2H(tp+(1-t)q)|dt.
\]
Using (\ref{udc2}), we estimate 
\begin{equation*}
\begin{split}
\langle D H(p)-D H(q), p-q\rangle &=~\int_0^1 \Big[ (p - q)^T A(t p+(1-t) q) (p - q) \Big]\cdot \big|D^2H (t p+(1-t) q)\big| dt \\[2mm]
%&~=~\int_0^1 \Big[ (p - q)^T A(t p+(1-t) q) (p - q) \Big]\big|D^2H (t p+(1-t) q)\big| dt \\[2mm]
&\geq~\lambda\cdot |p - q|^2 \int_0^1 \big|D^2H(t p+(1-t) q)\big| dt \\[2mm]
%&~\geq~{\lambda_1\over \lambda_2}\cdot |p - q| \left|\int_0^1  A( t p+(1-t) q) (p - q)  \big|D^2H (t p+(1-t) q)\big|  dt\right| \\[2mm]
&=~\lambda\cdot |D H(p)-D H(q)|\cdot |p -q|
\end{split}
\end{equation*}
and this implies (\ref{udc}).
\qed
\medskip

\noindent As we mentioned in the introduction, classical smooth solutions of~\eqref{HJ}  in general break down and Lipschitz continuous functions that satisfy~\eqref{HJ} almost everywhere together with a given initial condition  are not unique. To handle this problem, the following concept of a generalized solution 
was introduced in~\cite{CL}  to guarantee global existence and uniqueness results.
\begin{definition}\label{viscosity-solution}(Viscosity solution)
We say that a continuous function $u:[0,T]\times\mathbb{R}^d$ is a viscosity solution of \eqref{HJ} if:
\begin{enumerate}
\item[$\mathrm{(1)}$] u is a viscosity subsolution of \eqref{HJ}, i.e., for every point $(t_0,x_0)\in (0,T)\,\times\mathbb{R}^d$ and test function $v\in \mathcal{C}^1\big((0,+\infty)\times\mathbb{R}^d\big)$ such that $u-v$ has a local maximum at $(t_0,x_0)$, it holds that
\[
v_t(t_0,x_0)+H\big(D_{\!x} v(t_0,x_0)\big)\leqslant 0\,,
\] 
\item[$\mathrm{(2)}$]  u is a viscosity supersolution of (\ref{HJ}), i.e., for every point $(t_0,x_0)\in \,(0,T)\,\times\mathbb{R}^d$ and test function $v\in \mathcal{C}^1\big((0,+\infty)\times\mathbb{R}^d\big)$ such that $u-v$ has a local minimum at $(t_0,x_0)$, it holds that
\[
v_t(t_0,x_0)+H\big(D_{\!x} v(t_0,x_0)\big)\geqslant 0\,.
\] 
\end{enumerate}
\end{definition}
By the alternative equivalent definition of viscosity solution expressed in terms of the subdifferential and superdifferential of the function as in \cite{CL}  and because of Proposition \ref{con-c} one immediately  sees that every $\mathcal{C}^1$ solution of (\ref{HJ}) is also a viscosity solution of (\ref{HJ}). On the other hand, if $u$ is  a viscosity solution of (\ref{HJ}) then $u$ satisfies the equation at every point of differentiability. Let us state a result on further regularity for viscosity solutions proved in \cite[Proposition 3]{ACN} which says that  smoothness in the pair $(t,x)$ follows from smoothness in the second variable.
\begin{proposition}\label{conv-conx} Let $u$ be a viscosity solution of (\ref{HJ}) in $[0,T]\times\R^d$. If $u(t,\cdot)$ is both locally semiconcave and semiconvex in $\R^d$ for all $t\in (0,T)$ then $u$ is a $\mathcal{C}^1$ solution of (\ref{HJ}) in $(0,T]\times \R$.
\end{proposition}
The viscosity solution of the Hamilton-Jacobi equation (\ref{HJ}) with initial data $u(0,\cdot)=u_0\in \mathrm{{\bf Lip}}\big(\R^d\big)$
can be represented as the value function of a classical problem in calculus of variations, which admits the Hopf-Lax representation formula 
\begin{equation}\label{Hopf}
u(t,x)=\min_{y\in\mathbb{R}^d}\bigg\lbrace{t\cdot L\Big(\frac{x-y}{t}\Big)+u_0(y)\bigg\rbrace},
\qquad t>0, x\in \mathbb{R}^d,
\end{equation} 
where $L\in\mathcal{C}^1\left(\R^d\right)$ denotes the Legendre transform of $H$, defined by
\begin{equation}\label{Legendre-L}
L(q):=\max_{p\in\mathbb{R}^d}\big\{p\cdot q-H(p)\big\}, \qquad q\in\mathbb{R}^{d}.
\end{equation}
 The main properties of viscosity solutions defined by the Hopf-Lax formula, which are
of interest to this paper are recalled below (cfr. \cite[Section 1.1, Section 6.4]{CS}). 
\medskip
\begin{proposition}\label{ha-pr}
Let $u$ be the viscosity solution of \eqref{HJ}  on $[0,+\infty)\,\times\R^d$, with continuous initial data $u_0$,
defined by~\eqref{Hopf}. Then the following hold true:
\begin{enumerate}
\item[$(i)$] \textbf{Functional identity}: For all $x\in\mathbb{R}^d$ and $0\leqslant s<t$,
%\leqslant T$, 
it holds that 
\[
u(t,x)=\min_{y\in\mathbb{R}^d}\Big\lbrace{u(s,y)+(t-s)\cdot L\Big(\frac{x-y}{t-s}\Big)\Big\rbrace}\,.
\]
\item[$(ii)$] \textbf{Differentiability of $u$ and uniqueness}: \eqref{Hopf} admits a unique minimizer $y_x$ if and only if $u(t,\cdot)$ is differentiable at $x$. In this case we have 
\[
y_x=x-t\cdot DH\big(D_{\!x}u(t,x)\big),\qquad D_{\!x}u(t,x)\in D^-u_0(y_x).
\]
\item[$(iii)$] \textbf{Dynamic programming principle}: Let \, $t>s>0$, $x\in\mathbb{R}^d$,  assume that $y$ is a minimizer for \eqref{Hopf} and define $z=\ds \frac{s}{t}x+\left(1-\frac{s}{t}\right)y$. Then $y$ is the unique minimizer over~$\R^d$ of
\begin{equation*}
w\mapsto s\cdot L\Big(\frac{z-w}{s}\Big)+u_0(w),\qquad w\in\R^d.
\end{equation*}
\end{enumerate} 
\end{proposition}
\medskip

\noindent As a consequence, the family of nonlinear operators $\left\{S_t:{\bf Lip}\big(\mathbb{R}^d\big)\rightarrow {\bf Lip}(\mathbb{R}^d)\right\}_{t\geq 0}$
defined by the Hopf-Lax representation formula, i.e., $S_0 u_0=u_0$ and 
\bel{hopf-lax-smgr}
S_t u_0(x)=\min_{y\in\mathbb{R}^d}\Big\lbrace{t\cdot \ds L\Big(\frac{x-y}{t}\Big)+u_0(y)\Big\rbrace},\qquad
t>0\,,\;x\in\R^d,
\eeq
enjoys the following properties:
\begin{itemize}
\item[(i)] For every  $u_0\in{\bf Lip}\big(\mathbb{R}^d\big)$, $u(t,x):=S_t u_0(x)$ provides the unique viscosity solution of the 
Cauchy problem~\eqref{HJ} with initial data $u(0,\cdot)=u_0$.
\item[(ii)] (Semigroup property)
$$S_{t+s} u_0 = S_t\, S_s u_0,\qquad  t,s\geq 0, u_0\in{\bf Lip}(\R^d).$$
\item[(iii)] (Translation) For every constant $c\in\R$ we have that
\begin{equation}\label{eq:support}
S_t(u_0+c)=S_t u_0+c,\qquad t\geqslant 0,u_0\in{\bf Lip}(\R^d).
\end{equation}
%\item [(iv)] The map $S_t$ is continuous on sets of functions with uniform Lipschitz constant with respect to ${\bf W}^{1,1}_{\mathrm{loc}}$-topology, i.e., for every $u_n\in {\bf Lip}(\R^d)$ with a Lipschitz constant $M$ such that 
%\[
%u_n~\longrightarrow~u~~~\mathrm{in}~{\bf W}^{1,1}_{\mathrm{loc}}(\R^d),
%\]
%we have that  $S_t(u_n)$ also converges to $S_t(u)$ in ${\bf W}^{1,1}_{\mathrm{loc}}(\R^d)$.
\end{itemize}

\section{BV bound on ${\bf b}(t,\cdot )$}
Throughout this section, we shall assume that the Hamiltonian $H$ satisfies {\bf (H1)}-{\bf (H2)} and (\ref{Hs}). For a given initial datum $u_0\in {\bf Lip}\big(\R^d\big)$ with $\mathrm{Lip}[f]\leq M$, let $u$ be  the solution of (\ref{HJ}) with $u(0,\cdot)=u_0$ and
\bel{bd}
{\bf b}(t,x)=\left\{{x-y\over t}:y\in C_{t,x}\right\},\qquad C_{t,x}=\argmin_{y\in\R^d}\left\{u_0(y)+t\cdot L\left({x-y\over  t}\right)\right\}.
\eeq
It is well known from \cite{CS} that ${\bf b}(t,x)\subseteq DH(D^+u(t,x))$  and thus
\bel{LM}
\|{\bf b}(t,\cdot)\|_{{\bf L}^{\infty}(\R^d)}\leq\Lambda_{M}~:=~\max\{|q|:L(q)\leq M|q|\}.
\eeq
In order to establish a BV bound on  ${\bf b}(t,\cdot )$, we  approximate $u$ by a monotone decreasing sequence of continuous functions $u_n:[0,+\infty)\times \R^d\to\R$ defined by
\bel{u-n}
u_n(t,x):=\min_{y\in\mathcal{Z}_n} \left\{u_0(y)+t\cdot L\left({x-y\over t}\right)\right\},\qquad \mathcal{Z}_n:=2^{-n} \mathbb{Z}^d.
\eeq
Considering the associated set of slopes of backward optimal rays through $(t,x)$  
\bel{bn}
{\bf b}_n(t,x)=\left\{{x-y\over t}:y\in C^n_{t,x}\right\},\qquad C^n_{t,x}=\argmin_{y\in\mathcal{Z}_n}\left\{u_0(y)+t\cdot L\left({x-y\over  t}\right)\right\},
\eeq
we prove the following lemma.

\begin{lemma}\label{cv-un} For every $t>0$, it holds that  
\bel{bd-bb}
\lim_{n\to\infty}\|u_n(t,\cdot)-u(t,\cdot)\|_{\infty}=0\qquad\mathrm{and}\qquad \limsup_{n\to\infty}\|{\bf b}_n(t,\cdot)\|_{{\bf L}^{\infty}(\R^d)}~\leq~\Lambda_{M}.
\eeq
\end{lemma}
{\bf Proof.} {\bf 1.} Fix $n\geq 1$ and $(t,x)\in (0,\infty)\times\R^d$. Pick any $\bar{y}\in C_{t,x}$, let $\bar{z}$ be in $\mathcal{Z}_n$  such that  $|\bar{z}-\bar{y}|\leq \sqrt{d} 2^{-n+1}$. From (\ref{LM}) and (\ref{u-n}), we estimate 
\begin{eqnarray*}
|u(t,x)-u_n(t,x)|&\leq& \left|u_0(\bar{z})-u_0(\bar{y})\right|+t\cdot \left|L\left({x-\bar{z}\over  t}\right)-L\left({x-\bar{y}\over  t}\right)\right|\\
&\leq& \left(M+\max_{|q|\leq\Lambda_M+{\sqrt{d} \over 2^{n-1}\cdot t}} |DL(q)|\right)\cdot \big|\bar{z}-\bar{y}\big|
\end{eqnarray*}
and this yields the first equality in (\ref{bd-bb}).
\medskip

\n {\bf 2.} For  every $y_{x,n}\in C^n_{t,x}$ and $x'\in\R^d$, it holds that
\begin{eqnarray*}
u_n(t,x')-u_n(t,x)&\leq&t\cdot\left[L\left({x'-y_{x,n}\over t}\right)-L\left({x-y_{x,n}\over t}\right)\right]\\
&\leq&D L\left({x-y_{x,n}\over t}\right)\cdot (x'-x)+O(|x'-x|).
\end{eqnarray*}
In particular, ${\bf b}_n(t,x)\subseteq DH(D^+u_n(t,x))$ and the set  $\Sigma^n_{t}=\{x:\# {\bf b}_n(t,x)\geq 2\}$ is $\mathcal{H}^{n-1}$-rectifiable.  On the other hand, from (\ref{u-n})-(\ref{bn}), there exists $\bar{x}\in \mathcal{Z}_n $ with $|x-\bar{x}|\leq \sqrt{d} 2^{-n+1}$ such that  
\[
u_0(y_{x,n})+t\cdot L\left({x-y_{x,n}\over t}\right)\leq u_0(\bar{x})+t\cdot L\left({x-\bar{x}\over t}\right).
\]
Recalling that $\mathrm{Lip}[u_0]\leq M$, we get
\[
L\left({x-y_{x,n}\over t}\right)\leq M\cdot \left|{x-y_{x,n}\over t}\right|+\max_{|q|\leq {\sqrt{d} \over 2^{n-1}\cdot t}} L(q).
\]
From (\ref{Hs}), it holds that $L(0)=0$ and thus the above estimate yields the second part of (\ref{bd-bb}).
\qed
\medskip

\noindent Now we restate and prove Theorem \ref{BV-bound}, which is our first main theorem.
\begin{theorem} \label{newBV-bound} Assume that $H$ satisfies {\bf (H1)}-{\bf (H2)} and (\ref{Hs}).  For every $t>0$ and $u_0\in{\bf Lip}\big(\R^d\big)$ with $\mathrm{Lip}[u_0]\leq M$ for some $M\geq 0$, the function ${\bf b}(t,\cdot)$ has locally bounded total variation and its total variation $|D_x{\bf b}(t,\cdot)|$ over an open and bounded set $\Omega\subset\R^d$ of finite perimeter is bounded by 
\bel{BV-bb}
|D_x{\bf b}(t,\cdot)|(\Omega)~\leq~{1\over \gamma_{M}}\cdot \left(\Lambda_{M}+{\mathrm{diam}(\Omega)\over t}\right)\cdot\mathcal{H}^{d-1}(\partial\Omega)+{\sqrt{d}\over t}\cdot |\Omega|
\eeq
with $\gamma_M:=\ds\inf_{r>\max_{|q|\leq \Lambda_M}|DL(q)|}\lambda_r$ and $\lambda_r$ being as in (\ref{udc}).
\end{theorem}
\medskip

\n{\bf Proof.} The proof is divided into three steps\textcolor{red}{.}
\medskip

\n{\bf Step 1.} Consider the sequence of approximate solutions defined in (\ref{u-n}). Fixing  $n\geq 1$ and $t>0$, we write $\mathcal{Z}_n=\{y_1,y_2,\dots, y_k,\dots\}$. For any $i\neq j$, the set  
\[
\mathcal{O}_{i,j}=\left\{x\in\R^d:u_0(y_i)+t\cdot L\left({x-y_i\over t}\right)< u_0(y_j)+t\cdot L\left({x-y_j\over t}\right)\right\}
\]
is an open subset of $\R^d$ with a $\mathcal{C}^1$-boundary 
\[
\Gamma_{i,j}=\left\{x\in\R^d:u_0(y_i)+t\cdot L\left({x-y_i\over t}\right)=u_0(y_j)+t\cdot L\left({x-y_j\over t}\right)\right\}.
\]
Set $\mathcal{V}_i:=\bigcup_{j\neq i,j\geq 1} \mathcal{O}_{i,j}$. From  (\ref{u-n}) and (\ref{bn}), we have
\bel{dsb}
{\bf b}_n(t,x)={x-y_i\over t},\qquad x\in \mathcal{V}_i, ~i\geq 1.
\eeq
In particular, ${\bf b}_n(t,\cdot)$ is in $BV_{loc}(\R^d)$ and 
\bel{dbbb}
\begin{cases}
D_x{\bf b}_n(t,x)&=\ds{\mathbf{I}_d\over t}\cdot\mathcal{H}^d+{1\over 2t}\cdot \sum_{i\neq j} (y_j-y_i)\otimes \nu_i(x)\mathcal{H}^{d-1}_{\llcorner_{\partial{\mathcal{V}_i}\bigcap \partial{\mathcal{V}_j}}},\\
\div_x{\bf b}_n(t,x)&=~\ds{d\over t}\cdot \mathcal{H}^d+{1\over 2t}\cdot \sum_{i\neq j} \left\langle y_j-y_i, \nu_i(x)\right\rangle \mathcal{H}^{d-1}_{\llcorner_{\partial{\mathcal{V}_i}\bigcap \partial{\mathcal{V}_j}}},
\end{cases}
\eeq
where $\nu_i(x)$ is the inner normal vector to $\mathcal{V}_i$ and is computed by 
\[
\nu_i(x)~=~ {D L\big({x-y_j\over t}\big)-D L\left({x-y_i\over t}\right)\over \left|DL\big({x-y_j\over t}\big)-D L\big({x-y_i\over t}\big)\right|}\qquad\mathrm{for}~\mathcal{H}^{d-1}~a.e.~ x\in \partial{\mathcal{V}_i}\bigcap\partial{\mathcal{V}_j}.
\]
Given an open and bounded set $\Omega\subset\R^d$ of finite perimeter, one gets from (\ref{dbbb}) that
\bel{bDb}
\left|D_x{\bf b}_n(t,\cdot)-{\mathbf{I}_d\over t}\right|(\Omega)~\leq~{1\over 2t}\cdot \sum_{i\neq j} |y_j-y_i|\cdot\mathcal{H}^{d-1}(\Omega\bigcap\partial{\mathcal{V}_i}\bigcap \partial{\mathcal{V}_j}).
\eeq
For a fixed  $x\in \Omega\bigcap\partial{\mathcal{V}_i}\bigcap\partial{\mathcal{V}_j}$, setting $p_i:=D L\left({x-y_i\over t}\right)$ and $p_j:=D L\left({x-y_j\over t}\right)$, we have 
\[
\nu_i(x)={p_j-p_i\over |p_j-p_i|},\qquad y_j-y_i=D H(p_i)-D H(p_j).
\]
Recalling (\ref{udc}) and  (\ref{dsb}), we obtain that 
\bel{key}
|y_i-y_j|~\leq~-{1\over \lambda_{\beta_n}}\cdot \left\langle y_j-y_i, \nu_i(x)\right\rangle
\eeq
with $ {\beta_n}:=\ds\max_{|q|\leq \|{\bf b}_n(t,\cdot)\|_{{{\bf L}^{\infty}}(\R^d)}}|DL(q)|$ satisfying $\ds\lim_{n\to\infty}\beta_n=\max_{|q|\leq\Lambda_M}|DL(q)|$. Thus, (\ref{dbbb})-(\ref{bDb}) yield
\bel{est1}
\left|D_x{\bf b}_n(t,\cdot)-{\mathbf{I}_d\over t}\right|(\Omega)~\leq~{1\over \lambda_{\beta_n}}\cdot \left|\div_x{\bf b}_n(t,\cdot)-{d\over t}\right|(\Omega).
\eeq

\noindent {\bf Step 2.} Let us now provide a bound on $\left|\ds\div_x{\bf b}_n(t,\cdot)-{d\over t}\right|(\Omega)$, which will lead to a bound on $|D_x{\bf b}_n(t,\cdot)|(\Omega)$. Pick a point $x_0\in \Omega$. From (\ref{dsb}), (\ref{dbbb}) and (\ref{key}), the function ${\bf d}_n(t,\cdot):=\ds{\cdot-x_0\over t}-{\bf b}_n(t,\cdot)$ is in $BV_{\mathrm{loc}}\big(\R^d\big)$ and 
\[
\div_x {\bf d}_n (t,x)={1\over 2t}\cdot \sum_{i\neq j} \left\langle y_i-y_j, \nu_i(x)\right\rangle \mathcal{H}^{d-1}_{\llcorner_{\partial{\mathcal{V}_i}\bigcap \partial{\mathcal{V}_j}}}
\]
is a positive Radon measure. In particular, this implies that 
\[
|\div_x{\bf d}_n (t,\cdot)|(\Omega)=\int_{\Omega}\div_x{\bf d}_n (t,\cdot).
\]
Let $\rho_{\ve}\in\mathcal{C}^{\infty}_c\big(\R^d\big)$ be a family of modifiers, i.e., $\ds\rho_{\ve}(x)=\ve^{-d}\rho\left({x\over \ve}\right)$ for $\rho\in\mathcal{C}^{\infty}_c(\R^d)$ satisfying $\rho(x)\geq 0$, $\rho(x)=\rho(-x)$, $\mathrm{supp}(\rho)\subset B_d(0,1)$ and $\ds\int_{\R^d}\rho(x)dx=1$. For every  test function $\varphi_{\ve}=\chi_{\Omega}*\rho_{\ve}$, it holds that
\[
\int_{\R^d}\varphi_{\ve}\div_x{\bf d}_n (t,\cdot)=-\int_{\R^d} {\bf d}_n(t,x)\cdot \nabla\varphi_{\ve}(x)dx\leq\|{\bf d}_n(t,\cdot)\|_{{\bf L}^{\infty}(\R^d)}\cdot \int_{\R^d} | \nabla\varphi_{\ve}(x)|dx.
\]
Thus, taking $\ve\to 0+$, we get
\[
\int_{\Omega}\div_x{\bf d}_n (t,\cdot)~\leq~ \|{\bf d}_n(t,\cdot)\|_{{\bf L}^{\infty}(\R^d)}\cdot\mathcal{H}^{d-1}(\partial\Omega)
\]
and (\ref{est1}) yields
\bel{B-Db}
|D_x{\bf b}_n(t,\cdot)|(\Omega)~\leq~{1\over \lambda_{\Lambda_{M}}}\cdot \left(\|{\bf b}_n(t,\cdot)\|_{{\bf L}^{\infty}(\R^d)}+{\mathrm{diam}(\Omega)\over t}\right)\cdot\mathcal{H}^{d-1}(\partial\Omega)+{\sqrt{d}\over t}\cdot |\Omega|.
\eeq

\n {\bf Step 3.}  Finally, to achieve (\ref{BV-bb})  by taking $n\to\infty$ in (\ref{B-Db}), we first claim that ${\bf b}_n(t,\cdot)$ converges to ${\bf b}(t,\cdot)$ in ${\bf L}^{1}_{\mathrm{loc}}$. Since the sequence ${\bf b}_n(t,\cdot)$  is bounded in ${\bf L}^{\infty}\big(\R^d\big)$ and the set $\left(\bigcup_{n\geq 1}\Sigma^n_{t}\bigcup\Sigma_t\right)$ has zero Lebesgue measure with $\Sigma_t=\{x\in\R^d:\# {\bf b}(t,x)\geq 2\}$, it is sufficient to show that 
\[
\lim_{n\to\infty} {\bf b}_n(t,x)~=~{\bf b}(t,x)\qquad\forall x\in \R^d\backslash \left(\bigcup_{n\geq 1}\ds\Sigma^n_{t}\bigcup\Sigma_t\right).
\]
Given $x\in \R^d\backslash \left(\bigcup_{n\geq 1}\ds\Sigma^n_{t}\bigcup\Sigma_t\right)$, assume by a contradiction that there exists a subsequence ${\bf b}_{n_k}(t,x)$ converging to some $w\neq {\bf b}(t,x)$. From Lemma \ref{cv-un}, we have 
\begin{multline*}
u(t,x)~=~\lim_{n_k\to\infty} u_{n_k}(t,x)~=~\lim_{n_k\to\infty}u_0(x-t{\bf b}_{n_k}(t,x))+t\cdot L({\bf b}_{n_k}(t,x))\\
~=~u_0(x-tw)+t\cdot L(w)~=~u_0(x-tw)+t\cdot L\left({x- (x-tw)\over t}\right).
\end{multline*} 
Thus, ${\bf b}(t,x)$ is not a singleton and this yields a contradiction. By  \cite[Proposition 3.13]{ANP}, (\ref{B-Db}) and Lemma \ref{cv-un}, the function ${\bf b}_n(t,\cdot)$ converges weakly to ${\bf b}(t,\cdot)$ in $BV(\Omega,\R^d)$ with 
\begin{eqnarray*}
|D_x{\bf b}(t,\cdot)|(\Omega)&\leq&\liminf_{n\to\infty}|D_x{\bf b}_n(t,\cdot)|(\Omega)\\
&\leq&\left({1\over \ds\limsup_{n\to\infty}\lambda_{\beta_n}}\right)\cdot \left(\Lambda_M+{\mathrm{diam}(\Omega)\over t}\right)\cdot\mathcal{H}^{d-1}(\partial\Omega)+{\sqrt{d}\over t}\cdot |\Omega|
\end{eqnarray*}
and this yields (\ref{BV-bb}).
\qed 
\medskip

\noindent As a direct consequence of Theorem \ref{newBV-bound}, the following holds\textcolor{red}{.}
\medskip

\begin{corollary} The map $S_T:{\bf Lip}\left(\R^d\right)\to {\bf Lip}\left(\R^d\right)$ is compact in ${\bf W}^{1,1}_{{\mathrm{loc}}}\left(\R^d\right)$ for every  time $T>0$.
\end{corollary}
\section{Metric entropy in ${\bf W}^{1,1}$ for $S_T$}\label{metric-en}
In this section, we shall quantify  the degree of compactness of the map $S_T:{\bf Lip}\left(\R^d\right)\to {\bf Lip}\left(\R^d\right)$ for a given time $T>0$. More precisely, given constants $m,M,R>0$, considering the set of initial data
\bel{U-m}
\mathcal{U}_{[m,M]}=\left\{\bar{u}\in {\bf Lip}\big(\R^d\big): |\bar{u}(0)|\leq m, \mathrm{Lip}[\bar{u}]\leq M\right\},
\eeq
we establish upper and lower estimates for the $\ve$-entropy of the following restricted solution set at time $T$ in ${\bf W}^{1,1}(\square_R)$
\bel{S-TR}
S_{T,R}\left(\mathcal{U}_{[m,M]}\right):=\Big\{v_{\big\llcorner\square_R}~:~v\in S_T\left(\mathcal{U}_{[m,M]}\right)\Big\}.
\eeq
In order to do so, let us introduce continuous real-valued functions $\Psi_M,\Phi_M$ defined on $[0,M]$ for  $M>0$ such that $\Psi_M(0)=\Phi_M(0)=0$ and
\bel{Psi}
\begin{cases}
\Psi_M(s)&=\ds s\cdot \min_{p,q\in  \overline{B}_d(0,M), |p-q|\geq s} {|DH(p)-DH(q)|\over |p-q|}\\
\Phi_M(s)&=\ds s\cdot\min_{p\in \overline{B}_d\left(0,M-{s\over 2}\right)}\left(\max_{q\in \overline{B}_d(p,{s\over 2})} \left\|D^2H(q)\right\|_{\infty}\right)
\end{cases}
\quad\forall s\in (0,M].
\eeq
Here, $\left|D^2H(q)(v)\right|$ is the matrix norm of $D^2H(q)(v)$ and $\left\|D^2H(q)\right\|_{\infty}:=\ds\max_{|v|\leq 1}\left|D^2H(q)(v)\right|$. Notice that both maps $s\mapsto \Psi_M(s)$ and $s\mapsto \Phi_M(s)$ are strictly increasing. Moreover, the strict convexity   of $H$ implies that 
\[
0<\Psi_M(s)\leq\Phi_M(s)<M\cdot\max_{p\in \overline{B}_d(0,M)} \left\|D^2H(p)\right\|_{\infty}\quad\forall s\in (0,M].
\]

\noindent For convenience, we now rewrite Theorem \ref{second-main} as our second main theorem.

\begin{theorem}\label{main-2} Assume that $H\in\mathcal{C}^{2}\big(\R^d\big)$ and satisfies {\bf (H1)}-{\bf (H2)}. Then  for every $\ve>0$ sufficiently small, it holds that
\bel{sest1}
C_1\cdot \left(\ds \Phi_M\left({\ve\over R_1}\right)\right)^{-d}\leq\mathcal{H}_{\ve}\left(S_{T,R}(\mathcal{U}_{[m,M]})\Big|{\bf W}^{1,1}\big(\square_R\big)\right)~\leq~C_2\cdot \left(\ds \Psi_M\left({\ve\over R_2}\right)\right)^{-d}.
\eeq
for some constants  $C_1,C_2,R_1,R_2>0$ depending only on $m,M,R,T>0$.
\end{theorem}
\medskip

\noindent Before proving  Theorem \ref{main-2} in the next two subsections, we present some cases in which the estimates in (\ref{sest1}) are sharp.
\medskip

\begin{remark}  Given integer $k\geq 1$, the Hamiltonian $H(p)=|p|^{2k}$ is not uniformly convex but satisfies {\bf (H1)}-{\bf (H2)}. From (\ref{Psi}), there exist constants $\alpha_1,\alpha_2>0$ such that
\[
\alpha_1 s^{2k-1}\leq\Psi_M(s)\leq\Phi_M(s)\leq\alpha_2 s^{2k-1}\qquad\forall s\in [0,M].
\]
Thus, (\ref{sest1}) yields 
\[
\mathcal{H}_{\ve}\left(S_{T,R}(\mathcal{U}_{[m,M]})\big|{\bf W}^{1,1}\big(\square_R\big)\right)\approx\ve^{-(2k-1)d}.
\]
\end{remark}

\begin{remark} If $H\in\mathcal{C}^2(\R^d)$ is uniformly convex then $H$ satisfies {\bf (H1)}-{\bf (H2)} and
\[
\alpha_1s<\Psi_M(s)\leq\Phi_M(s)<\alpha_2 s \qquad\forall s\in [0,M],
\]
for some $\alpha_1,\alpha_2>0$. Thus,  (\ref{sest1}) yields the same result as in \cite{ACN} that 
\[
\mathcal{H}_{\ve}\left(S_{T,R}(\mathcal{U}_{[m,M]})\big|{\bf W}^{1,1}\big(\square_R\big)\right)\approx\ve^{-d}
\]
 for every $\ve>0$ sufficiently small.
\end{remark}

\noindent In the one-dimensional case, from Theorem \ref{main-2} we can obtain an estimate similar to that established in \cite[Remark 1.4]{AON4} for  scalar  conservation laws with strictly convex fluxes.

\begin{remark}
For $d=1$, every strictly convex $H\in\mathcal{C}^2(\R)$ satisfies (\ref{udc}). In addition, assume that $H$ has polynomial degeneracy, i.e., the set $I_H=\{\omega\in\R: H''(\omega)=0\}\neq \varnothing$ is finite and for each $w\in I_H$, there exists a natural number $p_\omega\geq 2$ such that 
\[
H^{(p_{\omega}+1)}(\omega)~\neq~0\qquad \mathrm{and}\qquad H^{(j)}(\omega)~=~0\quad\forall j\in \{2,\dots, p_{\omega}\}.
\]
The polynomial degeneracy of $H$  is defined by $\ds{\bf p}_{H}:=\max_{\omega\in I_H} p_{\omega}$.  For every $M>\ds\mathrm{argmax}_{w\in I_H}  p_{\omega}$, there exist constants $\alpha_1,\alpha_2>0$ such that 
\[
\alpha_1\cdot s^{{\bf p}_H}<\Psi_M(s)\leq\Phi_M(s)<\alpha_2\cdot s^{{\bf p}_H}\qquad\forall s\in [0,M].
\]
Thus, (\ref{sest1}) implies that  
$$\mathcal{H}_{\ve}\left(S_{T,R}(\mathcal{U}_{[m,M]})\big|{\bf W}^{1,1}\big(\square_R\big)\right)\approx \ve^{-{\bf p}_H}$$
 for every $\ve>0$ sufficiently small.
\end{remark}
\subsection{Upper estimate of  $\mathcal{H}_{\ve}\left(S_{T,R}(\mathcal{U}_{[m,M]})\big|{\bf W}^{1,1}\big(\square_R\big)\right)$}\label{UPP}
Towards the upper estimate of $\mathcal{H}_{\ve}\left(S_{T,R}(\mathcal{U}_{[m,M]})\big|{\bf W}^{1,1}\big(\square_R\big)\right)$ in (\ref{sest1}), we  first provide a bound on the ${\bf L}^1$-distance between elements $Du_1$ and $Du_2$ in terms of the ${\bf L}^1$-distance  between $DH(Du_1)$ and $DH(Du_2)$ for every $u_1,u_2\in S_{T,R}(\mathcal{U}_{[m,M]})$ by using the function $\Psi_M$ defined in (\ref{Psi}).  Observing that the map  $\ds s\mapsto { \Psi_M(s)\over s}$ is monotone increasing and 
\bel{Delta-p}
\Psi_M(|p-q|)\leq|DH(p)-DH(q)|\quad\forall p,q\in \overline{B}_d(0,M),
\eeq
we prove the following lemma.

\begin{lemma}\label{ds-L} For any $u_1,u_2\in S_{T,R}(\mathcal{U}_{[m,M]})$, it holds that
\bel{L1-et1}
\|Du_1-Du_2\|_{{\bf L}^{1}\big(\square_R\big)}\leq\left(2^dR^d+1\right)\cdot\Psi_M^{-1}\left(\|{\bf b}_1-{\bf b}_2\|_{{\bf L}^1\big(\square_R\big)}\right)
\eeq
with ${\bf b}_1:=DH(Du_1)$ and ${\bf b}_2:=DH(Du_2)$.
\end{lemma}
{\bf Proof.} For simplicity, setting $\alpha:=\Psi_M^{-1}\left(\|{\bf b}_1-{\bf b}_2\|_{{\bf L}^1\big(\square_R\big)}\right)$,
we claim that 
\bel{c1}
|Du_1(x)-Du_2(x)|\leq\alpha\cdot \max\left\{1, {|{\bf b}_1(x)-{\bf b}_2(x)|\over \|{\bf b}_1-{\bf b}_2\|_{{\bf L}^1\big(\square_R\big)}}\right\}\quad~\mathrm{for}~a.e.~ x\in \square_{R}.
\eeq
Indeed, assume that $|Du_1(x)-Du_2(x)|>\alpha$. From (\ref{Delta-p}), it holds that
\begin{multline*}
|Du_1(x)-Du_2(x)|={|Du_1(x)-Du_2(x)|\over |DH(Du_1(x))-DH(Du_2(x))|}\cdot |{\bf b}_1(x)-{\bf b}_2(x)|\\
\leq{|Du_1(x)-Du_2(x)|\over \Psi_M(|Du_1(x)-Du_2(x)|)}\cdot |{\bf b}_1(x)-{\bf b}_2(x)|.
\end{multline*}
By the monotone increasing property of the map $s\mapsto\ds \Psi_M(s)/s$, one has
\[
|Du_1(x)-Du_2(x)|\leq{\alpha\over |\Psi_M(\alpha)|}\cdot |{\bf b}_1(x)-{\bf b}_2(x)|~=~\alpha\cdot {|{\bf b}_1(x)-{\bf b}_2(x)|\over \|{\bf b}_1-{\bf b}_2\|_{{\bf L}^1\big(\square_R\big)}}
\]
and this implies (\ref{c1}). Hence, the ${\bf L}^1$-distance between $Du_1$ and $Du_2$ is bounded by 
 \begin{multline}\label{eee1}
\|Du_1-Du_2\|_{{\bf L}^1\big(\square_R\big)}=\int_{\square_R}|Du_1(x)-Du_2(x)|dx\\
\leq\alpha\cdot \int_{\square_R} \left(1+ {|{\bf b}_1(x)-{\bf b}_2(x)|\over \|{\bf b}_1-{\bf b}_2\|_{{\bf L}^1\big(\square_R\big)}}\right)dx=(2^dR^d+1)\alpha
\\ =\left(2^dR^d+1\right)\cdot\Psi_M^{-1}\left(\|{\bf b}_1-{\bf b}_2\|_{{\bf L}^1\big(\square_R\big)}\right)
\end{multline}
and this yields (\ref{L1-et1}).
\qed
\medskip

 {\bf  Proof of the upper estimate of  $\mathcal{H}_{\ve}\left(S_{T,R}(\mathcal{U}_{[m,M]})\big|{\bf W}^{1,1}\big(\square_R\big)\right)$ in  Theorem \ref{main-2}}
\v

\noindent {\bf 1.} From Theorem \ref{BV-bound}, for any $v\in S_{T,R}(\mathcal{U}_{[m,M]})$, one has  
\bel{bv-l-b}
|DH(Dv)|\left(\square_{R}\right)\leq V_T\qquad\mathrm{and}\qquad \|v\|_{{\bf L}^{\infty}(\R^d)}\leq m_T
\eeq
with  $V_T:=\ds{d2^d R^{d-1}\over \gamma_M}\cdot \left(\Lambda_{M}+{2\sqrt{d}R\over T}\right)+\ds{\sqrt{d}2^dR^d\over T}$ and $m_T:=\ds m+\sqrt{d}MR+T\cdot \sup_{|q|\leq \Lambda_M} L\left(q\right)$. In particular, the average value $\bar{v}^R$ of $v$ satisfies
\[
\bar{v}^R={1\over \mathrm{Vol}\big(\square_{R}\big)}\cdot \int_{\square_{R}} v(x)~dx~\in~[-m_T,m_T].
\]
Given $\ve'>0$, we  cover the interval  $[-m_T,m_T]$ by $K_{\ve'}=\ds\left\lfloor{m_T\over \Psi_M^{-1}(\ve')}\right\rfloor+1$ small intervals with length $2 \Psi_M^{-1}(\ve')$ such that 
\[
[-m_T,m_T]\subseteq\bigcup_{i=1}^{K_{\ve'}}B\big(a_i,\Psi_M^{-1}(\ve')\big)\quad\mathrm{for~some}~a_i\in [-m_T,m_T]
\]
and then decompose the set $S_{T,R}\big(\mathcal{U}_{[m,M]}\big)$  into $K_{\ve'}$ subsets as follows:
\[
S_{T,R}\left(\mathcal{U}_{[m,M]}\right)\subseteq\bigcup^{K_{\ve'}}_{i=1}\mathcal{A}_i,\qquad \mathcal{A}_i:=\left\{v\in S_{T,R}\left(\mathcal{U}_{[m,M]}\right): \bar{v}^R\in B\left(a_i,\Psi_M^{-1}(\ve')\right)\right\}.
\]
Recalling Definition \ref{DefKE}, we have 
\bel{N-e}
\mathcal{N}_{\ve}\left(S_{T,R}\left(\mathcal{U}_{[m,M]}\right)\Big|{\bf W}^{1,1}\big(\square_{R}\big)\right)\leq\sum_{i=1}^{K_{\ve'}}\mathcal{N}_{\ve}\left(\mathcal{A}_i\Big|{\bf W}^{1,1}\big(\square_{R}\big)\right)
\eeq
for all $\ve>0$.
\v

\noindent {\bf 2.} For a given $i\in\{1,2,\dots,K_{\ve'}\}$, we are going to  provide an upper bound on the covering number $\mathcal{N}_{\ve}\left(\mathcal{A}_i\Big|{\bf W}^{1,1}\big(\square_{R}\big)\right)$
by introducing the  set $\mathcal{B}_i:=\left\{DH(Dv): v\in \mathcal{A}_i\right\}$. From (\ref{F}) and (\ref{bv-l-b}), one has that $\mathcal{B}_i\subseteq\mathcal{F}_{\left[R,\Lambda_M,V_{T}\right]}$ with $\mathcal{F}_{\left[R,\Lambda_M,V_{T}\right]}$ defined as in (\ref{F}). By  Corollary \ref{e-p}, if $\ve'>0$ is sufficiently small then it holds that
\bel{H-e}
\mathcal{H}_{\ve'/2}\left(\mathcal{B}_i\big|{\bf L}^1\left(\square_R\right)\right)\leq\Gamma^+\cdot(\ve')^{-d},\qquad \Gamma^+=48\sqrt{d}\cdot \left(12d\sqrt{d}RV_T\right)^d.
\eeq
By the definition of $\mathcal{H}_{\ve'/2}\left(\mathcal{B}_i\big|{\bf L}^1\left(\square_R\right)\right)$, there exists a set $\left\{{\bf v}_1,\dots, {\bf v}_{\beta_{\ve'}}\right\}\subset\mathcal{A}_i$ with $\beta_{\ve'}\leq \ds2^{\Gamma^+\cdot (\ve')^{-d}}$ such that \[
\mathcal{B}_i\subseteq\bigcup_{j=1}^{\beta_{\ve'}} B_{{\bf L}^1}({\bf b}_j,\ve'),\qquad {\bf b}_j~:=~DH(D{\bf v}_j).
\]
In particular, for any given $v\in \mathcal{A}_i$, it holds that
\[
\|DH(Dv)-{\bf b}_{j_0}\|_{{\bf L}^1\big(\square_{R}\big)}~<~\ve'~~~\mathrm{for~some}~j_0\in \overline{1,\beta_{\ve'}}.
\]
Recalling Lemma \ref{ds-L}, we obtain that
\begin{eqnarray*}
\|Dv-D{\bf v}_{j_0}\|_{{\bf L}^1\big(\square_R\big)}&\leq&\left(2^dR^d+1\right)\cdot\Psi_M^{-1}\left(\|DH(Dv)-{\bf b}_{j_0}\|_{{\bf L}^1\big(\square_R\big)}\right)\\
&\leq&\left(2^dR^d+1\right)\cdot\Psi_M^{-1}(\ve')
\end{eqnarray*}
and the Poincar\'e inequality in Theorem \ref{PC} yields
\begin{eqnarray*}
\left\|\left(v-\bar{v}^R\right)- \left({\bf v}_{j_0}-\bar{{\bf v}}_{j_0}^R\right)\right\|_{{\bf L}^1\big(\square_R\big)}&\leq&\sqrt{d}R\cdot \|Dv-D{\bf v}_{j_0}\|_{{\bf L}^1\big(\square_R\big)}\\
&\leq&\sqrt{d}R\left(2^dR^d+1\right)\cdot\Psi_M^{-1}(\ve').
\end{eqnarray*}
On the other hand, since $v, {\bf v}_{j_0}\in \mathcal{A}_i$, one has 
\[
\big|\bar{v}^R-\bar{{\bf v}}_{j_0}^R\big|\leq\big|\bar{v}^R-a_{i}\big|+\big|\bar{{\bf v}}_{j_0}^R-a_{i}\big|~\leq~2\Psi_M^{-1}(\ve').
\]
Thus, the ${\bf W}^{1,1}$-distance between $v$ and ${\bf v}_{j_0}$ can be estimated by  
\begin{eqnarray*}
\|v-{\bf v}_{j_0}\|_{{\bf W}^{1,1}\big(\square_R\big)}&\leq& \big|\bar{v}^R-\bar{{\bf v}}_{j_0}^R\big|\cdot \big|\square_R\big|+\|Dv-D{\bf v}_{j_0}\|_{{\bf L}^1\big(\square_R\big)}\\
&+&\left\|\left(v-\bar{v}^R\right)- \left({\bf v}_{j_0}-\bar{{\bf v}}_{j_0}^R\right)\right\|_{{\bf L}^1\big(\square_R\big)}\leq R^+\cdot \Psi_M^{-1}(\ve')
\end{eqnarray*}
with $R^+:=\left(2^dR^d +1\right)(3+\sqrt{d}R)$. Finally, by choosing $\ve' = \ds\Psi_M\left({\ve\over R^+}\right)$, we have that  $\ds\mathcal{A}_i\subseteq\bigcup_{i=1}^{\beta_{\ve'}}B_{{\bf W}^{1,1}}({\bf v}_i,\ve)$ and 
\[
\mathcal{N}_{\ve}\left(\mathcal{A}_i\big|{\bf W}^{1,1}\big(\square_{R}\big)\right)\leq\beta_{\ve'}~=~\ds 2^{\ds\Gamma^+\cdot \left(\Psi_M\left({\ve\over R^+ }\right)\right)^{-d}}.
\]
Therefore, from (\ref{N-e}), one gets
\[
\mathcal{N}_{\ve}\left(S_{T,R}(\mathcal{U}_{[m,M]})\Big|{\bf W}^{1,1}\big(\square_{R}\big)\right)\leq\ds\left(\left\lfloor{m_TR^+\over \ve}\right\rfloor+1\right)\cdot 2^{\ds\Gamma^+\cdot \left({ \Psi_M\left({\ve\over R^+}\right)}\right)^{-d}}
\]
and this yields the second inequality in (\ref{sest1}) for $\ve>0$ sufficiently small.
\qed
\begin{remark} To obtain the upper bound of $\mathcal{H}_{\ve}\left(S_{T,R}(\mathcal{U}_{[m,M]})\big|{\bf W}^{1,1}\big(\square_R\big)\right)$ in (\ref{sest1}), we only require that $H$ belongs to $\mathcal{C}^1(\R^d)$ and satisfies {\bf (H1)-(H2)}. 
\end{remark}
\subsection{Lower estimate of  $\mathcal{H}_{\ve}\left(S_{T,R}(\mathcal{U}_{[m,M]})\big|{\bf W}^{1,1}\big(\square_R\big)\right)$}\label{LUP} In this subsection, we shall prove the first inequality in (\ref{sest1}). In order to do so, for any given  $p\in\R^d$, let  $\Phi(\cdot, p):[0,\infty) \to [0,\infty)$ be the strictly increasing continuous function defined by  $\Phi(0,p)=0$ and 
\[
\Phi(s,p)=s\cdot \left(\max_{p'\in \overline{B}_d(p,{s\over 2})} \left\|D^2H(p')\right\|_{\infty}\right),\qquad s>0.
\]
From the definition of $\Phi_M$ in (\ref{Psi}), it holds that
\bel{cdn}
\Phi_{M}(s)=\min_{p\in\overline{B}_d\left(0,M-{s\over 2}\right)}\Phi(s,p),\qquad s\in [0,M].
\eeq
Let us recall the constant in the assumption {\bf (H2)}
\bel{udc1}
\lambda_r=\inf_{p\neq q\in\overline{B}(0,r)}\left\langle {DH(p)-DH(q)\over |DH(p)-DH(q)|}, {p-q\over |p-q|}\right\rangle>0,\qquad r>0.
\eeq
The following proposition shows that a solution to (\ref{HJ})  with a semiconvex initial condition preserves the semiconvexity on a given time interval,
provided the semiconvexity constant of the initial data is sufficiently small in absolute value. 

\begin{proposition}\label{r-convex} Assume that $H\in\mathcal{C}^{2}\big(\R^d\big)$ and satisfies {\bf (H1)}-{\bf (H2)}. Given $T,M,r>0$ and $\bar{p}\in\overline{B}_d\left(0,M-{r\over 2}\right)$, let $\bar{u}$ be a semiconvex function with semiconvexity constant $- K$ such that 
\bel{u0-c}
D^-\bar{u}(\R^d)\subseteq\overline{B}_d\left(\bar{p},{r\over 2}\right),\qquad K\leq{\lambda_M\over 4T}\cdot {r\over \Phi(r,\bar{p})}.
\eeq
Then, the map $(t,x)\mapsto S_t(\bar{u})(x)$ is a classical solution for $0<t\leq T$ and 
\[
DS_t(\bar{u})(x)\in\overline{B}_d\left(\bar{p},{r\over 2}\right)\qquad\forall (t,x)\in (0,T]\times \R^d.
\]
\end{proposition}

\noindent {\bf Proof.} For  simplicity, we set 
\[
u(t,x):=S_t(\bar{u})(x)\qquad\forall (t,x)\in [0,\infty)\times\R^d.
\]
It is well-known from \cite[Theorem 5.3.8]{CS} that $u(t,\cdot)$ is locally semiconcave for every $t>0$. Thus, by Proposition \ref{conv-conx}, it is sufficient to show that 
$u(t,\cdot)$ is  semiconvex with some semiconvexity constant $-C<0$ for all $t\in [0,T]$, i.e., for any fixed $(t,x)\in [0,T)\times \R^d$, it holds that
\bel{le1}
u(t,x+h)+u(t,x-h)-2u(t,x)\geq -C\cdot |h|^{2}\quad\forall h\in \R^d.
\eeq
By the Lipschitz continuity of $u(t,\cdot)$, we can assume that $u(t,\cdot)$ is differentiable at $x\pm h$.  In this case, ${\bf b}(t,x\pm h)$ reduce to a single value denoted by ${\bf b}^{\pm}=DH({\bf p}^{\pm})$ with ${\bf p}^{\pm}=Du(t,x\pm h)$ and satisfy the following relations
\begin{equation}\label{pp}
\begin{cases}
{\bf p}^{\pm}\in D^-\bar{u}(x \pm h-t {\bf b}^{\pm})\subseteq\ds\overline{B}_d\left(\bar{p},{r\over 2}\right)\subseteq\overline{B}_d(0,M),\cr\cr
u(t,x\pm h)=\bar{u}(x \pm h-t{\bf b}^{\pm})+ t\cdot L({\bf b}^{\pm}).
\end{cases}
\end{equation}
Since $\bar{u}$ is semiconvex with semiconvexity constant $-K$, denoting $x^{\pm}:=x\pm h$, from (iii) of Proposition \ref{con-c} one can get that
\[
\big\langle {\bf p}^{+}-{\bf p}^{-}, x^+-x^-- t\left({\bf b}^+-{\bf b}^-\right)\big\rangle \geq-K\cdot \left|2h-t\left({\bf b}^+-{\bf b}^-\right)\right|^2
\]
and 
\begin{eqnarray*}
\big\langle  {\bf p}^{+}-{\bf p}^{-},{\bf b}^+-{\bf b}^-\big\rangle&\leq&{K\over t}\cdot \left|2h-t\left({\bf b}^+-{\bf b}^-\right)\right|^2+{2|h|\over t}\cdot |{\bf p}^{+}-{\bf p}^{-}|\\
&\leq& 2Kt|{\bf b}^{+}-{\bf b}^{-}|^2+{8K|h|^2\over t}+{2|h|\over t}\cdot |{\bf p}^{+}-{\bf p}^{-}|\\
&\leq&  2KT |DH({\bf p}^+)-DH({\bf p}^-)|^2+{8K|h|^2\over t}+{2|h|\over t}\cdot |{\bf p}^{+}-{\bf p}^{-}|.
\end{eqnarray*}
Since ${\bf p^{\pm}}\in \overline{B}_d(\bar{p},{r\over 2})$, it holds that
\[
|DH({\bf p}^+)-DH({\bf p}^-)|\leq{\Phi(r,\bar{p})\over r}\cdot |{\bf p}^+-{\bf p}^-|.
\]
Thus, recalling  (\ref{u0-c}) and (\ref{pp}), we estimate 
\begin{multline*}
 2KT |DH({\bf p}^+)-DH({\bf p}^-)|^2\leq2KT\cdot {\Phi(r,\bar{p})\over r}\cdot  |DH({\bf p}^+)-DH({\bf p}^-)|\cdot |{\bf p}^+-{\bf p}^-|\\
 \leq{\lambda_M\over 2}\cdot |DH({\bf p}^+)-DH({\bf p}^-)|\cdot |{\bf p}^+-{\bf p}^-|={\lambda_M\over 2}\cdot |{\bf b}^+-{\bf b}^{-}|\cdot |{\bf p}^+-{\bf p}^{-}|
\end{multline*}
and 
\bel{one-side}
\big\langle  {\bf p}^{+}-{\bf p}^{-},{\bf b}^+-{\bf b}^-\big\rangle \leq{\lambda_M\over 2}\cdot  |{\bf b}^+-{\bf b}^{-}|\cdot |{\bf p}^+-{\bf p}^-|+{8K|h|^2\over t}+{2|h|\over t}\cdot |{\bf p}^{+}-{\bf p}^{-}|.
\eeq
On the other hand, from (\ref{udc1}) we deduce  that 
\begin{multline*}
\big\langle  {\bf p}^{+}-{\bf p}^{-},{\bf b}^+-{\bf b}^-\big\rangle=\big\langle  {\bf p}^{+}-{\bf p}^{-},DH({\bf p}^+)-DH({\bf p}^-)\big\rangle\\
\geq\lambda_M\cdot |DH({\bf p}^+)-DH({\bf p}^-)|\cdot |{\bf p}^{+}-{\bf p}^{-}|=\lambda_M\cdot |{\bf b}^+-{\bf b}^{-}|\cdot |{\bf p}^+-{\bf p}^{-}|
\end{multline*}
and  (\ref{one-side}) implies that 
\bel{ke1}
{\lambda_M\over 2}\cdot  |t ({\bf b}^+-{\bf b}^{-})|\cdot |{\bf p}^+-{\bf p}^{-}|~\leq~8K|h|^2+2|h| \cdot |{\bf p}^{+}-{\bf p}^{-}|.
\eeq
Observe that if $|{\bf p}^+-{\bf p}^{-}|\leq K|h|$ then   
\begin{eqnarray*}
|{\bf b}^+-{\bf b}^{-}|= |DH({\bf p}^+)-DH({\bf p}^-)|\leq{\Phi(r,\bar{p})\over r}\cdot |{\bf p}^+-{\bf p}^-|\leq {K\Phi(r,\bar{p})\over r}\cdot |h|.
\end{eqnarray*}
Otherwise, (\ref{ke1}) implies that $\ds {\lambda_M\over 2}\cdot |t ({\bf b}^+-{\bf b}^{-})|\leq10|h|.$ Hence, it holds that
\bel{ke2}
 |t ({\bf b}^+-{\bf b}^{-})|~\leq~\left({KT\Phi(r,\bar{p})\over r}+{20\over \lambda_M}\right)\cdot |h|.
\eeq
By the Hopf-Lax representation formula, we have 
\[
u(t,x\pm h)~=~\bar{u}(x\pm h-t{\bf b}^{\pm})+t\cdot L({\bf b}^{\pm})
\]
\[
u(t,x)~\leq~\-2\bar{u}\left(x-t\cdot{{\bf b}^++{\bf b}^-\over 2}\right)+t\cdot L\left({\bf b}^++{\bf b}^-\over 2\right).
\]
Using the convexity of $L$ and semiconvexity of $\bar{u}$, we estimate 
\begin{multline*}
u(t,x+h)+u(t,x-h)-2u(t,x)\geq t\cdot \left[L({\bf b}^+)+L({\bf b}^-)-2L\left({\bf b}^++{\bf b}^-\over 2\right)\right]\\
+\bar{u}(x+h-t{\bf b}^+)+\bar{u}(x-h-t{\bf b}^-)-2\bar{u}\left(x-t\cdot{{\bf b}^++{\bf b}^-\over 2}\right)\\
\geq-K\cdot \left|2h-t\left({\bf b}^+-{\bf b}^-\right)\right|^2\geq-8K|h|^2-2K |t({\bf b}^+-{\bf b}^{-})|^2\\
\geq- 2K\cdot \left[4+\left({KT\Phi(r,\bar{p})\over r}+{20\over \lambda_M}\right)^2\right]\cdot |h|^2
\end{multline*}
and this yields  (\ref{le1}).
\qed
\medskip

\noindent Relying on the above Proposition and Corollary \ref{low}, we now proceed to prove the first inequality in (\ref{sest1}).
\v

 {\bf  Proof of the lower estimate of $\mathcal{H}_{\ve}\left(S_{T,R}(\mathcal{U}_{[m,M]})\big|{\bf W}^{1,1}\big(\square_R\big)\right)$ in  Theorem \ref{main-2}}
\v

\noindent {\bf 1.} Let us recall that $\mathcal{U}_{[m,M]}=\left\{\bar{u}\in {\bf Lip}\big(\R^d\big): |\bar{u}(0)|\leq m, \mathrm{Lip}[\bar{u}]\leq M\right\}$ and 
\[
\mathcal{SC}_{[r,K]}:=\left\{v\in {\bf Lip}(\R^d): \mathrm{Lip}[v]\leq r~\mathrm{and}~ v~\mathrm{is~semiconcave~with~constant}~K\right\}.
\]
For any given $r>0$ and $p\in \R^d$, we  denote by 
\[
\mathcal{W}^{p}_{r}:=\left\{  \varphi= v+\langle p,\cdot\rangle:v\in \mathcal{SC}_{[{r\over 2}, K_r]}\right\},\qquad K_r={\lambda_M\over 4T}\cdot {r\over \Phi_M(r)}.
\]
 From  (\ref{cdn}), there exists $p_r\in \overline{B}_d(0,M-{r\over 2})$ such that $\Phi(r,p_r)=\ds\min_{p\in\overline{B}\left(0,M-{r\over 2}\right)}\Phi(r,p)$. Now consider the operator $\mathcal{T}:\mathcal{W}^{p_r}_{r}\to {\bf Lip}\big(\R^d\big)$ such that for all $\varphi\in \mathcal{W}^{p}_{r}$, it holds that
\bel{T}
\mathcal{T}(\varphi)~=~\varphi+S_T(\varphi_-)(0),\qquad  \varphi_-(\cdot):=-\varphi(-\cdot).
\eeq
By reversing the equation (\ref{HJ}), we will show that   
\bel{inc22}
\mathcal{T}\left(\mathcal{W}^{p_r}_{r}\right)\subseteq S_T\left(\mathcal{U}_{[0,M]}\right).
\eeq
Indeed, for a given $\varphi\in\mathcal{W}^{p_r}_{r}$, we define the following function  
\[
w_0(\cdot):=-\mathcal{T}(\varphi)(-\cdot)=\varphi_-(\cdot)-S_T(\varphi_-)(0).
\]
Since $\varnothing\neq D^+\varphi(x)\subset\overline{B}_d \left(0,{r\over 2}\right)$ for all $x\in\R^d$, $w_0$  is semiconvex with a semiconvexity constant $-K_r$ and 
\[		
D^-w_0(x)=p_r+D^+\varphi(-x)\subseteq\overline{B}_d \left(p_r,{r\over 2}\right)\qquad\forall x\in \R^d.
\]
Let $w(t,x)=S_t(w_0)(x)$ be the unique viscosity solution of  (\ref{HJ}) with initial datum $w(0,x)=w_0$. Recalling Proposition \ref{r-convex} and property (ii) in Proposition \ref{ha-pr}, we have that $w$ is a $\mathcal{C}^1$ classical solution of (\ref{HJ}) in $(0,T]\times \R^d$ and 
\[
D_xw(T,x)\subseteq\overline{B}_d\left( p_r,{r\over 2}\right)\subseteq\overline{B}_d(0,M)\qquad\forall x\in \R^d.
\]
Moreover, from the translation property of $S_t$ in  (\ref{eq:support}), it holds that
\[
w(T,0)=S_T(w_0)(0)=S_T\left(\varphi_--S_T(\varphi_-)(0)\right)(0)=0.
\]
Thus, the continuous function $u:[0,T]\times\R^d\to\R$, defined by 
\[
u(t,x)=-w(T-t,-x)\qquad\forall (t,x)\in [0,T]\times\R^d,
\]
is also a $\mathcal{C}^1$ classical solution of (\ref{HJ}) in $(0,T)\times\R^d$  with
\[
u(T,\cdot)=\mathcal{T}(\varphi)(\cdot)\quad\mathrm{and}\quad u(0,\cdot)=-w(T,-\cdot)\in\mathcal{U}_{[0,M]}.
\] 
In particular, $u(t,x)$ is  a viscosity solution of (\ref{HJ}) in $[0,T]\times\R^d$, so that by the
uniqueness property of the semigroup map $S_t$, one has 
\[
S_T(u_0)(\cdot)=\mathcal{T}(\varphi)(\cdot),\qquad u_0(\cdot)=-w(T,-\cdot)
\]
and this yields (\ref{inc22}).
\v

\noindent {\bf 2.} For every $\ve>0$, we select a finite subset $A_{\ve}\subseteq [-m,m]$ such that 
\bel{ca1}
\# A_{\ve}~=~\left\lfloor{2^dR^d m\over \ve}\right\rfloor\quad\mathrm{and}\quad |a_i-a_j|~\geq~{2\ve\over 2^dR^d}\quad\forall a_i\neq a_j \in A_{\ve}.
\eeq
Again from the translation property of $S_t$ in  (\ref{eq:support}), one has
\[
S_T\big(\mathcal{U}_{[m,M]}\big)\supseteq\bigcup_{a\in A_{\ve}} S_T\big(a+\mathcal{U}_{[0,M]}\big)~=~A_{\ve}+S_T\big(\mathcal{U}_{[0,M]}\big)
\]
and  (\ref{inc22})  implies that 
\bel{inc3}
S_T\big(\mathcal{U}_{[m,M]}\big)\supseteq A_{\ve}+\mathcal{T}\left(\mathcal{W}^{p_r}_{r}\right).
\eeq
By Corollary \ref{low},  for every $\ve>0$ sufficiently small, there exists a set $\mathcal{G}\subseteq\mathcal{W}^{p_r}_{r}$ such that 
\[
\# \mathcal{G}\geq\ds 2^{\beta_1\cdot \ve^{-d}},\qquad \beta_{1}={1\over 3^{d}2^{d^2+4d+3}\ln 2}\cdot \left(\omega_d R^{d+1}K_r\over(d+1)\right)^d
\]
and 
\[
\left\|D\varphi_{\big\llcorner\square_R}-D\phi_{\big\llcorner\square_R}\right\|_{{\bf L}^1\big(\square_R\big)}\geq2\ve\qquad\forall \varphi\neq\phi\in \mathcal{G}.
\]
Since $D\mathcal{T}(\varphi)(x)=D\varphi(x)$ for all $x\in\R^d$, 
\[
\left\|D\mathcal{T}(\varphi)_{\big\llcorner\square_R}-D\mathcal{T}(\phi)_{\big\llcorner\square_R}\right\|_{{\bf W}^{1,1}\big(\square_R\big)}\geq 2\ve\qquad\forall \varphi\neq\phi\in \mathcal{G}.
\]
Recalling (\ref{ca1}), we  have 
\[
\left\|f_{\big\llcorner\square_R}-g_{\big\llcorner\square_R}\right\|_{{\bf W}^{1,1}\big(\square_R\big)}\geq2\ve\quad\forall f\neq g\in A_{\ve}+\mathcal{T}\left(\mathcal{W}^{p_r}_{r}\right)
\]
and thus
\begin{eqnarray*}
\mathcal{H}_{\ve}\left(S_{T,R}(\mathcal{U}_{[m,M]})\Big|{\bf W}^{1,1}\big(\square_R\big)\right)\geq\log_2\left(\# A_{\ve}\cdot \# \mathcal{G}\right)=\log_2\left(\left\lfloor{2^dR^d m\over \ve}\right\rfloor\right)+{\beta_1\over \ve^d}.
\end{eqnarray*}
Finally, by choosing $\ds r={\ve\over R^-}$ with $R^-:=\ds{\omega_d\cdot R^d\over (d+1)2^{d+9}}$, we compute 
\[
K_r={\lambda_M\over 4T R^-}\cdot {\ve\over \Phi_M\left({\ve\over R^-}\right)},\qquad \beta_{1}~=~{1\over 8\ln2}\cdot \left(8R\lambda_M\over 3T\right)^d\cdot \left({\ve\over \Phi_M\left({\ve\over R^-}\right)}\right)^d
\]
and get
\[
\mathcal{H}_{\ve}\left(S_{T,R}\big(\mathcal{U}_{[m,M]}\big)\Big|{\bf W}^{1,1}\big(\square_R\big)\right)\geq {1\over 8\ln2}\cdot \left(8R\lambda_M\over 3T\right)^d \cdot \left( \Phi_M\left({\ve\over R^-}\right)\right)^{-d}+\log_2\left(\left\lfloor{2^dR^d m\over \ve}\right\rfloor\right).
\]
This particularly yields the first inequality in (\ref{sest1}) for every $\ve>0$ sufficiently small.
\qed
\begin{remark} With the same argument,  the lower bound of  $\mathcal{H}_{\ve}\left(S_{T,R}(\mathcal{U}_{[m,M]})\big|{\bf W}^{1,1}\big(\square_R\big)\right)$ in (\ref{sest1}) can be obtained for $H\in\mathcal{C}^{1,1}(\R^d)$ by defining 
\[
\Phi_M(s)=s\cdot\inf_{p\in \overline{B}_d\left(0,M-{s\over 2}\right)}\left(\sup_{p_1\neq p_2\in \overline{B}_d(p,{s\over 2})} {|DH(p_1)-DH(p_2)|\over |p_1-p_2|}\right)
\] 
for all $s>0$.
\end{remark}

\section{A counter-example}
\label{Ss:No BV}
In this section, we  provide an example to show that if  the  Hamiltonian $H\in\mathcal{C}^2(\R^2)$ satisfies {\bf (H1)} but not {\bf (H2)} then Theorem \ref{BV-bound} fails in general. Consider a smooth, coercive and strictly convex Hamiltonian
\[
H(p)={3^3\over 4^4}\cdot p_1^4+p_2^2,\qquad p=(p_1,p_2)\in\R^2.
\]
The function $H$  does not satisfy {\bf (H2)} as
$$\ds\lim_{p_1\to 0}\left\langle {DH(p_1,p^2_1)-DH(0,0)\over |DH(p_1,p_1^2)-DH(0,0)|}, {(p_1,p^2_1)\over |(p_1,p^2_1)|}\right\rangle= \ds\lim_{p_1\to 0}\left\langle {\left({3^3\over 4^3}\cdot p_1^3,2p_1^2\right)\over \left|\left({3^3\over 4^3}\cdot p_1^3,2p_1^2\right)\right|}, {(p_1,p^2_1)\over |(p_1,p^2_1)|}\right\rangle=0.$$
In the next steps, we shall construct an  initial datum $u_0\in\mathrm{{\bf Lip}}(\R^2)$ such that both $D_{\!x}u(1,\cdot)$ and ${\bf b}(1,\cdot)=DH(D_{\!x}u(1,\cdot))$  do not have locally bounded  variation where $u$ is the solution of (\ref{HJ}) with $u(0,\cdot)= u_0$.
\v

\n{\bf Step 1:} For given $0<\ell<1$,  we first construct an initial datum $\bar{u}\in\mathrm{{\bf Lip}}(\R^2)$  with  $\mathrm{Lip}[\bar{u}]\leq 1$  such that 
\begin{equation}\label{Equa:bd1}
\mathrm{supp}(\bar{u})\subset[-2\ell,2\ell],\qquad |{\bf b}(1,\cdot)|([-\ell,\ell]^2), |D_{\!x}u(1,\cdot)|([-\ell,\ell]^2)\geq1
\end{equation}
where   $u$ is  the solution of (\ref{HJ}) with $u(0,\cdot)= \bar{u}$. For every $0<\delta<\ell$, we consider the periodic lattice 
\[
y_{\iota}=\left(\iota_1 \delta, \iota_2\delta^{2/3}\right),\qquad \iota\in \mathcal{Z}_2:=\left\{(\iota'_1,\iota'_2)\in \mathbb{Z}^2: \iota'_1+\iota'_2\in 2\mathbb{Z}\right\}
\]
and the  corresponding regions 
\begin{eqnarray*}
\Omega_{\iota}&=&\{x\in \R^2:L(x-y_{\iota})< L(x-y_{\iota'})\quad\forall \iota'\neq \iota\}\\
&=&y_{\iota}+\left\{q\in\R^2: L(q)< L(q+ y_{\iota}-y_{\iota'})\quad\forall \iota'\neq \iota\right\}\\
&\subseteq& y_{\iota}+ [-\delta,\delta]\times [-\delta^{2/3},\delta^{-2/3}].
\end{eqnarray*}
Let $g_0,g_1:\R^2\to\R$ be such that 
\[
\begin{cases}
g_1(x)&=L(x-y_{\iota}),\qquad x\in \Omega_{\iota}, \iota\in\mathcal{Z}_2,\\
g_0(y)&=~\max\limits_{x\in\R^2} \left\{g_1(x)-L(x-y)\right\},\quad y\in\mathbb{R}^2.
\end{cases}
\]
Notice that  both $g_0$ and $g_1$ are  Lipschitz  with Lipschitz constant 
\[
M_{\delta}~=~\sup_{q\in [-\delta,\delta]\times [-\delta^{2/3},\delta^{-2/3}]} |DL(q)|~=~O(\delta^{1/3}).
\]
Thus, for  $\delta>0$  sufficiently small, one can construct $\bar{u}\in\mathrm{{\bf Lip}}(\R^2)$ with  $\mathrm{Lip}[\bar{u}]\leq M_{\delta}$ and 
\[
\mathrm{supp}(\bar{u})~\subset~[-2\ell,2\ell],\qquad \bar{u}(y)~=~g_0(y)~~\forall y\in \left[-{3\ell\over 2},{3\ell\over 2}\right]^2.
\]
Let $u$  be  the solution of (\ref{HJ}) with $u(0,\cdot)= \bar{u}$. At time $t=1$, we have  
\[
u(1,x)~=~\min_{y\in\R^2} \left\{\bar{u}(y)+L(x-y)\right\}~=~\bar{u}(y_x)+ L(x-y_x)
\]
for some $y_x\in \overline{B}(x,\Lambda_{M_{\delta}})$ with $\Lambda_{M_{\delta}}=\max\{|q|: L(q)\leq M_{\delta}\cdot |q|\}=O(\delta^{1/3})$. Thus, if $M_{\delta}\leq \ds{\ell\over 2}$ then for all $x\in [-\ell,\ell]^2\cap \Omega_{\iota}$, $\iota\in\mathcal{Z}_2$,
\begin{eqnarray*}
u(1,x)&=&\min_{y\in \left[-{3\ell\over 2},{3\ell\over 2}\right]^2} \left\{\bar{u}(y)+L(x-y)\right\}~=~\min_{y\in \left[-{3\ell\over 2},{3\ell\over 2}\right]^2} \left\{g_0(y)+L(x-y)\right\}\\
&=&\min_{y\in\R^2} \left\{g_0(y)+L(x-y)\right\}~=~g_1(x)~=~L(x-y_{\iota})
\end{eqnarray*}
and the slope of  backward optimal rays through $(1,x)$ is
\[
{\bf b}(1,x)~=~DH(D_{\!x}u(1,x))~=~x-y_{\iota}.
\]
 For any two adjacent $y_{\iota},y_{\iota'}$ with $\Omega_{\iota},\Omega_{\iota'}\subset [-\ell,\ell]^2$ and $x\in\partial\Omega_{\iota}\cap \partial\Omega_{\iota'}$, denoting the inner normal vector to $\Omega_\iota$ by ${\bf n}(x)$, we compute
\begin{equation*}
\begin{cases}
D_{\!x}u(1,x)&=~[DL(x-y_{\iota})-DL(x-y_{\iota'})]\otimes {\bf n}(x) \mathcal H^1_{\llcorner_{\partial \Omega_\iota \cap \partial \Omega_{\iota'}}},\cr\cr
D_{\!x}{\bf b}(1,x)&=~ (y_{\iota}-y_{\iota'})\otimes {\bf n}(x) \mathcal H^1_{\llcorner_{\partial \Omega_\iota \cap \partial \Omega_{\iota'}}}.
\end{cases}
\end{equation*}
From the definition of $\Omega_{\iota}$, one can show that $\mathcal{H}^1(\partial\Omega_{\iota}\cap \partial\Omega_{\iota'})\geq \delta^{2/3}$ and this implies 
\begin{equation*}
|D_{\!x}{\bf b}(1,\cdot)|(\Omega_{\iota}\cup\Omega_{\iota'}), |D_{\!x}u(1,\cdot)|(\Omega_{\iota}\cup\Omega_{\iota'})~\geq~  \delta^{2/3}\cdot \mathcal{H}^1(\partial\Omega_{\iota}\cap \partial\Omega_{\iota'})~\geq~ \delta^{4/3}.
\end{equation*}
Moreover, since the number of open regions  $\Omega_{\iota}\subset [-\ell,\ell]^2$ is of the order $\ds {\ell^2\over \delta^{5/2}}$, there exists a constant $C>0$ such that
\begin{equation*}
|D_{\!x}{\bf b}(1,\cdot)|([-\ell,\ell]^2), |D_{\!x}u(1,\cdot)|([-\ell,\ell]^2)~\geq~C\cdot {\ell^2\over \delta^{5/2}}\cdot \delta^{4/3}~=~C\cdot{\ell^2\over \delta^{1/3}}.
\end{equation*}
Thus, choosing $\delta>0$ sufficiently small, we obtain \eqref{Equa:bd1}.
\v

\noindent {\bf Step 2.} Consider a sequence of disjoint squares  $\square_{n}=c_n+[0,2^{-n}]\times [0,2^{-n}]$ such that $\ds \bigcup_{n\geq 1}\square_{n}\subset [0,1]^2$.
From the previous step, for any $n\geq 1$ one can construct a sequence of functions  $u_{0,n}\in \mathrm{{\bf Lip}}(\R^2)$  with $\mathrm{Lip}[u_{0,n}]\leq 1$ such that $\mathrm{supp}(u_{0,n})\subset\square_n$ and the solution $u_n$ of \eqref{HJ} with  $u_n(0,\cdot)=u_{0,n}(\cdot)$ satisfies  
\[
|D_{\!x} u_n(1,\cdot)|\left(c_n+{1\over 2}\cdot\left(\square_n-c_n\right)\right),~|DH(D_{\!x}u_n(1,\cdot))|\left(c_n+{1\over 2}\cdot\left(\square_n-c_n\right)\right)~\geq~1
\]
and 
\[
L(x-z)~\geq~\min_{y\in \square_n}\left\{u_{0,n}(y)+L(x-y)\right\},\quad x\in\left(c_n+{1\over 2}\cdot\left(\square_n-c_n\right)\right), z\in\R^2\backslash\square_n\,.
\]
Finally, set $u_0=\ds\sum_{n=1}^{\infty}u_{0,n}\in \mathrm{{\bf Lip}}(\R^2)$. The solution $u$ of \eqref{HJ} with $u(0,\cdot)=u_0(\cdot)$ satisfies 
\[
u(1,x)~=~u_n(1,x),\qquad x\in \left(c_n+{1\over 2}\cdot\left(\square_n-c_n\right)\right)
\]
and this implies  
\[
|D_{\!x}{\bf b}(1,\cdot )|([0,1]^2)~\geq~\sum_{n=1}^{\infty}|DH(D_{\!x} u_n(1,\cdot))|(\square_n)~\geq~\sum_{n=1}^{\infty} 1~=~+\infty.
\]
Similarly, one has that $|D_{\!x}u(1,\cdot)|([0,1]^2)=+\infty$.
\qed
\v

{\bf Acknowledgments.} This research by K.T. Nguyen was partially supported by a grant from
the Simons Foundation/SFARI (521811, NTK). The authors would like to warmly thank the anonymous referees for carefully reading the manuscript and for their suggestions, which greatly helped in improving the paper overall.


\begin{thebibliography}{99}

\bibitem{AD}G. Acosta and R. C. D\'uran, An optimal Poincar\'e inequality in ${\bf L}^1$ for convex domains,
{\it Proc. Amer. Math. Soc. Vol} {\bf 132} (2003), no.1, p. 195-202.

\bibitem{ANP} L. Ambrosio, N. Fusco and D. Pallara, Functions of Bounded Variation and Free Discontinuity Problems, {\it Oxford Science Publications}, Clarendon Press, Oxford, UK, (2000).
\
\bibitem{ACN} F. Ancona, P. Cannarsa and K. T. Nguyen, Quantitative compactness estimates for Hamilton-Jacobi equations, {\it Arch. Rat. Mech. Anal.}, {\bf 219} (2016), no. 2, 793--828. 
\bibitem {ACN1} F. Ancona, P. Cannarsa and K. T. Nguyen, The compactness estimates for Hamilton Jacobi Equations depending on space, {\it Bull. Inst. Math. Acad. Sin. (N.S.)} {\bf 11} (2016), no. 1, 63--113.
\bibitem{AON1} F. Ancona, O. Glass and K. T. Nguyen, Lower compactness estimates for scalar balance laws, {\it Comm. Pure Appl. Math} {\bf 65} (2012), no. 9, 1303--1329.
\bibitem{AON3} F. Ancona, O. Glass and K. T. Nguyen, On quantitative compactness estimates for hyperbolic conservation laws,
{\it Proceedings of the 14th International Conference on Hyperbolic Problems (HYP2012)}, AIMS, Springfield, MO, (2014).
\bibitem{AON2} F. Ancona, O. Glass and K. T. Nguyen, On lower compactness estimates for general nonlinear hyperbolic systems, {\it Ann. Inst. H. Poincare Anal. Non Lineaire}, {\bf 32} (2015), no. 6, 1229--1257. 

\bibitem{AON4}  F. Ancona, O. Glass and K. T. Nguyen, On Kolmogorov entropy compactness estimates for scalar conservation laws without uniform convexity, {\it SIAM J. Math. Anal.}, {\bf 51} (2019), no. 4, 3020--3051.

 
\bibitem{CS} P. Cannarsa and C. Sinestrari, Semiconcave functions,
Hamilton-Jacobi equations, and optimal control,  {\it Progress in Nonlinear
Differential Equations and their Applications}, 58. Birkh\"auser
Boston, (2004). 
%

 \bibitem{DNR}  R. Capuani, P. Dutta and K. T. Nguyen, Metric entropy for functions of bounded total generalized variation, {\it SIAM J. Math. Anal.} {\bf 53} (2021), no. 1, 1168--1190.

\bibitem{CL} M.G. Crandall and P.-L. Lions, Viscosity solutions of Hamilton-Jacobi equations, {\it Trans. Amer. Math. Soc.} {\bf 277} (1983), no. 1, 1--42.

 
\bibitem{DLG} C. De Lellis and F. Golse, A Quantitative Compactness Estimate for Scalar Conservation Laws, {\it Comm. Pure Appl. Math.} {\bf 58} (2005), no. 7, 989--998. 

 \bibitem{DN} P. Dutta and K. T. Nguyen,  Covering numbers for bounded variation functions, {\it J. Math. Anal. Appl.} {\bf 468} (2018), no. 2, 1131--1143.

\bibitem{KT} A.N. Kolmogorov and V.M Tikhomirov, $\varepsilon$-Entropy and $\varepsilon$-capacity of sets in functional spaces,
{\it Uspekhi Mat. Nauk} {\bf 14} (1959), 3-86.

\bibitem{Lax02} P.D. Lax, Course on hyperbolic systems of conservation laws, {\it XXVII Scuola Estiva di Fis. Mat.}, Ravello, 2002.

\end{thebibliography}
\end{document}